\documentclass[11pt]{article}
\usepackage{amsmath}
\usepackage{amssymb}
\usepackage{bbm}
\usepackage[mathscr]{euscript}
\usepackage[margin=2cm]{geometry}
\usepackage{hyperref}
\usepackage{tikz}
\usetikzlibrary{decorations.pathmorphing,arrows}
\usepackage{braket}
\usepackage{setspace}

\newtheorem{thm}{Theorem}[section]

\newtheorem{lemma}[thm]{Lemma}

\newtheorem{proposition}[thm]{Proposition}

\newtheorem{corollary}[thm]{Corollary}

\newenvironment{proofof}{}{\hfill$\square$\vskip.5cm}
\newcommand{\R}{\mathbb{R}}
\newcommand{\N}{\mathbb{N}}
\newcommand{\C}{\mathbb{C}}

\newcommand{\Z}{\mathbb{Z}

}

\renewcommand{\Re}{\operatorname{Re}}
\renewcommand{\Im}{\operatorname{Im}}
\newcommand{\x}{\boldsymbol{x}}
\newcommand{\y}{\boldsymbol{y}}

\renewcommand{\Z}{{\mathbb{Z}}}

\begin{document}

\onehalfspace

%
\title{About the Moments of the Generalized Ulam Problem}
\date{May 9, 2024}


\author{
Samen Hossein\\
Department of Economics\\
University of Pennsylvania\\
133 South 36th Street
Suite 150\\
Philadelphia, PA 19104
\\[7pt]
Shannon Starr\\
 University of Alabama at Birmingham\\
 1402 Tenth Avenue South\\
 Birmingham, AL 35294-1241\\
 \href{mailto:slstarr@uab.edu}{slstarr@uab.edu}}

\maketitle

\abstract{Given $\pi \in S_n$, let  $Z_{n,k}(\pi)=\sum_{1\leq i_1<\dots<i_k\leq n} \mathbf{1}(\{ \pi_{i_1}<\dots<\pi_{i_k}\}$ denote the number of 
increasing subsequences of length $k$. 
Consider the``generalized Ulam problem,'' studying the distribution of  $Z_{n,k}$ for general $k$ and $n$.
For the 2nd moment, Ross Pinsky initiated a combinatorial study by considering a pair of subsequences $i^{(r)}_1<\dots<i^{(r)}_k$ for $r \in \{1,2\}$,
and conditioning on the size of the intersection $j = |\{i_1^{(1)},\dots,i^{(1)}_k\} \cap \{i^{(2)}_1,\dots,i^{(2)}_k\}|$.
We obtain the exact large deviation rate function for $\mathbf{E}[Z_{n,k} Z_{n,\ell}]$ in the asymptotic regime $k\sim \kappa n^{1/2}$, $\ell \sim  \lambda n^{1/2}$ as $n \to \infty$,
for $\kappa,\lambda \in (0,\infty)$. This uses multivariate generating function techniques, as found in the textbook  of Pemantle and Wilson.
The requisite generating function enumerates pairs of up-right paths in $d=2$, which both end at $(k,\ell)$ with a given number of intersections.
We also evaluate the analogous generating function for pairs of $(+\boldsymbol{i},+\boldsymbol{j},+\boldsymbol{k})$ paths in $d=3$, which both end at $(k,\ell,m)$,
which has some utility in calculating the 3rd moment.

Finally, we consider a simpler problem involving partitions instead of permutations, where all moments are calculable and the replica symmetric ansatz can be stated if not proved.
}

\thispagestyle{empty}

\section{Introduction}



Suppose we have an integer $n \in \N = \{1,2,\dots\}$ and another integer $k \in [n] = \{1,\dots,n\}$.
Let $S_n$ denote the set of all permutations of $\pi = (\pi_1,\dots,\pi_n)$ (i.e., such that $\{\pi_1,\dots,\pi_n\} = [n]$). Let us denote a function $Z_{n,k} : S_n \to \R$ defined as 
\begin{equation}
	Z_{n,k}(\pi)\, =\, \sum_{i_1=1}^{n-k+1}  \sum_{i_2=i_1+1}^{n-k+2} \cdots \sum_{i_k=i_{k-1}+1}^{n}  \mathbf{1}(\{\pi_{i_1}<\dots<\pi_{i_k}\})\, .
\end{equation}
In other words, $Z_{n,k}(\pi)$ is equal to the number of subsequences of $\pi$ which are of length $k$ and which are increasing relative to $\pi$.
We call the ``generalized Ulam problem'' the investigation of the joint distribution of the random variables $Z_{n,k}(\pi)$ for $k \in [n]$,
when $\pi \in S_n$ is uniform random.

In \cite{LifschitzPittel}, Lifschitz and Pittel considered the distribution of $X_n=\sum_{k=1}^{n} Z_{n,k}$.
In his textbook \cite{Steele}, Steele uses the $Z_{n,k}$'s to show that the length of the longest increasing subsequence 
of $\pi$ concentrates around a number $\ell_n$ which is asymptotically $c n^{1/2}$ for some $c \in [1,e]$.
The actual $c$ is $2$, as determined by Vershik and Kerov \cite{VershikKerov} and independently Logan and Shepp \cite{LoganShepp} (whose results were more partial).
The fluctuations of the length of the longest increasing subsequence were determined by Baik, Deift and Johansson \cite{BDJ}
using the RSK algorithm and a Riemann-Hilbert problem.
It is the same as the fluctuations of the largest eigenvalue of a GUE random matrix: the Tracy-Widom distribution \cite{TracyWidom}.
The investigation of the length of the longest increasing subsequence is known as the Ulam problem, since it was initiated by him, with later contributions by Hammersley.
Talagrand's method gives a concentration of measure bound for the length of the longest increasing subsequence quite easily \cite{Talagrand}.
But it was surprisingly recently, with the work of Baik, Deift and Johansson, that the problem of determining the fluctuations was solved.

The reason for calling the investigation of the $Z_{n,k}$'s the ``generalized Ulam problem'' is that the length of the longest increasing subsequence of $\pi$ is 
$\max(\{k \in [n]\, :\, Z_{n,k}(\pi)>0\})$. In other words, if $L_n(\pi)$ is the length of the longest increasing subsequence of $\pi$, then 
we have $L_n(\pi) \geq k$ if and only if $Z_{n,k}(\pi) > 0$.
In a pair of important papers, Ross Pinsky considered when $Z_{n,k}$ satisfies a weak law of large numbers.
In his first paper, he gave a sufficient condition for $Z_{n,k}$ to satisfy a WLLN, using the Payley-Zygmund second moment method \cite{Pinsky1}.
Then, using the Baik-Deift-Johansson result, in a second paper he gave a sufficient condition for $Z_{n,k}$ to fail to satisfy a WLLN.
The two conditions are not quite matching, but they do establish that $Z_{n,k}$ fails to satisfy a WLLN whenever $k = c n^{1/2}$ for any $c$ (which is the most interesting
case from the perspective of Ulam's original problem).

Pinsky did not calculate the precise asymptotics of the 2nd moment of $Z_{n,k}$, although he did obtain rough bounds, suitable for his purpose.
Moreover, he did obtain a direct combinatorial approach to the variance.
Here we generalize his combinatorial approach, which is suitable for calculating covariances between $Z_{n,k}$ and $Z_{n,\ell}$ for $k\neq \ell$.

We also initiate the problem of computing higher moments for the generalized Ulam problem. We propose it as an interesting for future analysis.
The reason for considering higher moments is that this is the first step in the physicists' replica method.
In order to demonstrate this on a similar problem, we introduce a truly solvable model in section \ref{sec:solvable}.

One might think in rough analogy to the Sherrington-Kirkpatrick model, which turned out to be solvable but only with a great extra effort.
In comparison, the Random Energy Model is immediately solvable, because the energy levels are uncorrelated.
The solvable model we introduce is akin to a version of the generalized Ulam problem where the correlations are much simpler 
(just as the correlations in the REM are much simpler than the SK model).

We discuss what we can say for the 3rd moment of the generalized Ulam problem in Section \ref{sec:higher}.

\section{Main results for asymptotics of the 2nd moments}
\subsection{Second moments : combinatorics}
The main combinatorial result is as follows.
\begin{proposition}
\label{prop:PinskyGeneralization}
For any $n \in \N = \{1,2,\dots\}$ and any $k \in [n]$, it is known $\mathbf{E}[Z_{n,k}] = \binom{n}{k} / k!$.
Now, for any $k,\ell \in [n]$, 
\begin{equation}
	\mathbf{E}[Z_{n,k}Z_{n,\ell}]\, =\, \sum_{j=0}^{k\wedge \ell} \mathbf{E}[Z_{n,k+\ell-j}] \mathcal{A}(k-j,\ell-j,j)\, ,
\end{equation}
where 
\begin{equation}
	\mathcal{A}(k,\ell,j)\, =\, \sum_{\substack{\alpha_0,\dots,\alpha_j \in \{0,1,\dots\}\\ \alpha_0+\dots+\alpha_j=k}} 
 \sum_{\substack{\beta_0,\dots,\beta_j \in \{0,1,\dots\}\\ \beta_0+\dots+\beta_j=\ell}} \prod_{r=0}^j \binom{\alpha_r+\beta_r}{\alpha_r,\beta_r}^2\, . 
\end{equation} 
\end{proposition}
This proposition is a minor generalization of a result of  Pinsky, Proposition 2 in his article \cite{Pinsky1}.

The way to think of this is as follows. Sum over all  pairs of subsequences $1\leq a_1<\dots<a_k \leq n$ and $1\leq b_1<\dots<b_\ell\leq n$.
Given such a pair, let $j = |\{a_1,\dots,a_k\}\cap \{b_1,\dots,b_\ell\}|$. 
Then let $1\leq c_1<\dots<c_{k+\ell-j}\leq n$ be the sequence such that
$\{c_1,\dots,c_{k+\ell-j}\} = \{a_1,\dots,a_k\}\cap \{b_1,\dots,b_\ell\}$. Let us condition on $j$ and $(c_1,\dots,c_{k+\ell-j})$,
to turn around the procedure we just described. Careful consideration of the combinatorics involved in such a move leads to the proposition. A full proof will be
provided in Section \ref{sec:Pictures}.

Pinsky also showed a direct bijection to obtain a follow-up formulation.
We state the following version.

\begin{corollary}
\label{cor:RW}
Suppose that $(W_t)_{t=0}^{k+\ell}$ is a 2d simple random walk, such that $W_0=(0,0)$ and $W_{t}-W_{t-1}$ are IID uniform random variables in the set 
$\{(1,0) ,(-1,0),(0,1),(0,-1)\}$. Then
\begin{equation}
	\mathcal{A}(k,\ell,j)\, =\, 4^{k+\ell} \sum_{0\leq t_0\leq \dots \leq t_j\leq k+\ell}
 \mathbf{P}\Big(W_{k+\ell} = (k-\ell,0)\ \text{ and }\
W_{t_0},\dots, W_{t_j} \in \Z \times \{0\}\Big)\, .
\end{equation}
\end{corollary}

The bijection that Pinsky gives is a bijection between two 1d independent SRW's $(X_0,\dots,X_{k+\ell})$ and $(Y_0,\dots,Y_{k+\ell})$
and a 2d SRW $(W_0,\dots,W_{k+\ell})$ obtained from the mapping $W_t^{(1)} = (X_t+Y_t)/2$ and $W_t^{(2)}=(X_t-Y_t)/2$. 
Note for instance that then the probability that $W_{k+\ell} = (k-\ell,0)$ is equal to $P(X_{k+\ell}=Y_{k+\ell}=k-\ell) = 4^{-(k+\ell)} \binom{k+\ell}{k,\ell}^2$.
(In mathematical physics this might be considered a diagonal transfer or 45-degree rotation of the usual orientation as in Baxter \cite{Baxter}.)
We will give a proof of this corollary in Section \ref{sec:Pictures}.

Using the diagonal transfer one can also prove the generating function formula: for $x,y \in \mathcal{U}(0;1/4)$:
\begin{equation}
\label{eq:CombGF}
	\sum_{\alpha=0}^\infty \sum_{\beta = 0}^{\infty} \binom{\alpha+\beta}{\alpha,\beta}^2 
x^{\alpha} y^{\beta}\, =\, \frac{1}{\big(1 - 2(x+y) + (x-y)^2\big)^{1/2}}\, .
\end{equation}
This uses the diagonal transfer method and the Newton version of the binomial formula for the square-root function.
Then by a geometric series type calculation one obtains the following.
\begin{corollary}
\label{cor:GF}
For $x,y,z \in \mathcal{U}(0;r)$ for $r=\sqrt{5}-2$ (where $\mathcal{U}(\zeta_0;\rho) = \{\zeta \in \C\, :\, |\zeta-\zeta_0|<\rho\}$) we have
\begin{equation}
	\sum_{k=0}^{\infty} \sum_{\ell=0}^{\infty} \sum_{j=0}^{\infty} \mathcal{A}(k,\ell,j) x^k y^{\ell} z^j\, =\, \frac{1}{\big(1 - 2(x+y) + (x-y)^2\big)^{1/2} - z}\, .
\end{equation}
\end{corollary}

We will present a proof in Section \ref{sec:Pictures}.

\subsection{Asymptotics of the combinatorial factor in the critical regime}

In this section, we state the asymptotics in their simplest form of large deviation rate functions.
We will describe corrections in a later section.

The critical regime is $k \sim \kappa n^{1/2}$ as $n \to \infty$ for $\kappa \in (0,\infty)$.
The reason is that then 
\begin{equation}
\label{eq:FirstMomentAsymp}
\lim_{n \to \infty} n^{-1/2} \ln(\mathbf{E}[Z_{n,k}])\, =\, 2 \kappa \big(1-\ln(\kappa)\big)\, .
\end{equation}
because some larger terms cancel.
One can see that this becomes negative for $\kappa>e$ meaning the expected value of $Z_{n,k}$ is less than 1 for $k/n^{1/2} \sim \kappa > e$.
In particular that suggests the length of the longest increasing subsequence is no longer than $e n^{1/2}$ asymptotically as $n \to \infty$,
and that can be made rigorous as in Steele \cite{Steele}.

Then it is known that 
\begin{equation}
\label{eq:SymmetricAAsymptotics}
\begin{split}
	\left(\lim_{K \to \infty} \frac{j_K}{4K}\, =\, \rho \in (0,\infty)\right) \qquad \Rightarrow\\
&\hspace{-3cm} \qquad \left(\lim_{K \to \infty} \frac{1}{2K}\, \ln\left(\mathcal{A}(K,K,j_K)\right)\, =\, 2 \ln(2) + (1+\rho) \ln(1+\rho) - \rho \ln(\rho)\right)
\end{split}
\end{equation}
This is an exact version of a rough bound from \cite{Pinsky1}.
However, using Corollary \ref{cor:GF}, one can derive the following generalization.

\begin{thm}
Given sequences $k_N$, $\ell_N$ and $j_N$, assuming $\lim_{N \to \infty} k_N/N = \kappa$, $\lim_{N \to \infty} \ell_N/N = \lambda$ and $\lim_{N \to \infty} j_N/N = \gamma$
for $\kappa,\lambda,\gamma \in (0,\infty)$,  we have
\begin{equation}
	\lim_{N \to \infty} N^{-1} \ln(\mathcal{A}(k_N,\ell_N,j_N))\, =\, -\kappa \ln(X) - \lambda \ln(Y) - \gamma \ln(Z)\, ,
\end{equation}
for $X = X(\kappa,\lambda,\gamma)$, $Y = Y(\kappa,\lambda,\gamma)$ and $Z = Z(\kappa,\lambda,\gamma)$ as below:
\begin{equation}
X\, =\, \frac{\kappa \cdot (2\kappa+\gamma)}{(\kappa+\lambda+\gamma)(2\kappa+2\lambda+\gamma)}\, ,\qquad
Y\, =\, \frac{\lambda \cdot (2\lambda+\gamma)}{(\kappa+\lambda+\gamma)(2\kappa+2\lambda+\gamma)}\, ,
\end{equation}
and
\begin{equation}
Z\, =\, \frac{1}{\kappa+\lambda+\gamma}\, \cdot 
\left(\frac{\gamma\cdot (2\kappa+\gamma)(2\lambda+\gamma)}{2\kappa+2\lambda+\gamma}\right)^{1/2}\, .
\end{equation}
\end{thm}
In other words, for $\psi(x) = -x \ln(x)$, and $h(\theta) = \psi(\theta) + \psi(1-\theta)$, we have
\begin{equation}
\label{eq:LastRewriteAAsymptotics}
\begin{split}
	(\kappa+\lambda)^{-1} \lim_{N \to \infty}  N^{-1} \ln(\mathcal{A}(\kappa N,\lambda N,2\rho(\kappa+\lambda) N))\, 
&=\, 2\, \ln(2) - \psi(1+\rho) + \psi(\rho) \\
&\qquad+ h\left(\frac{\kappa}{\kappa+\lambda}\right) - h(1/2) \\
&\hspace{-3cm}
+(1+2\rho) \Big( h\left(\frac{1+\rho}{1+2\rho} \cdot \frac{\kappa}{\kappa+\lambda} + \frac{\rho}{1+2\rho} \cdot \frac{\lambda}{\kappa+\lambda}\right) - h(1/2)\Big)\, ,
\end{split}
\end{equation}
where $h(1/2) = \ln(2)$ is the maximum value of $h(\theta)$ for $\theta \in [0,1]$.
The 2nd and 3rd lines on the right-hand-side of equation (\ref{eq:LastRewriteAAsymptotics})
evaluate to $0$ if $\kappa=\lambda$. Therefore, those lines can be viewed as corrections to equation (\ref{eq:SymmetricAAsymptotics})
for the case $\kappa \neq \lambda$.

Note that with the formulas given $(X,Y,Z)$ lies on the affine variety $z^2 = 1+2(x+y)+(x-y)^2$ which is the singularity set of the generating function
$\sum_{j,k,\ell} \mathcal{A}(k,\ell,j) x^k y^\ell z^j$. 
One step in the determination of which point on the variety to choose for the optimal $(X,Y,Z)$ is to use coordinates that are as homogeneous as possible.

Once we have chosen the optimizer $(X,Y,Z)$, our basic tool is the saddle point method. In the univariate case, a classic reference is the textbook of Flajolet and Sedgewick \cite{FlajoletSedgewick}.
An important reference for multivariate generating functions is Pemantle and Wilson \cite{PemantleWilson}.

\subsection{Asymptotics of the second moment in the critical regime}

For simplicity, we continue to state the results in the large deviation form (without describing the algebraic factors that correct the basic exponential behavior).

\begin{lemma}
If we have a sequence $k_n$ satisfying $k_n \sim \kappa n^{1/2}$ as $n \to \infty$, for some $\kappa \in (0,\infty)$, then
\begin{equation}
	\frac{1}{2\kappa} \lim_{n \to \infty}  n^{-1/2} \ln\left(\mathbf{E}\left[Z_{n,k}^2\right] \right)\,
=\, 1-2\ln(2)+\frac{1}{1+4P}-\frac{1}{2}\, \ln(P\cdot (1+P))\, ,
\end{equation}
where $P=P(\kappa) \in (0,\infty)$ is the unique solution of the equation 
\begin{equation}
\label{eq:r2kappaP}
\frac{2P^{1/4}(1+P)^{3/4}(1+4P)}{1+2P}\, =\, \kappa\, .
\end{equation}
\end{lemma}
This is a straightforward application of the Laplace method or Varadhan's lemma using equation \ref{eq:SymmetricAAsymptotics}.
See, for example, \cite{DemboZeitouni}.
The value of $P(\kappa)$ gives the asymptotics of the $j$ which maximizes the summand in the formula:
\begin{equation}
	\mathbf{E}[Z_{n,k}^2]\, =\, \sum_{j=0}^{k} \mathbf{E}[Z_{n,2k-j}] \mathcal{A}(k-j,k-j,j)\, .
\end{equation}
Namely, if we have
$k \sim \kappa n^{1/2}$ as $n \to \infty$, then the optimal $j$ is such that $j \sim 4 P(\kappa) \kappa n^{1/2}$ as $n \to \infty$.


Note that Lifschitz and Pittel in \cite{LifschitzPittel} proved that the logarithm of 
the 2nd moment
for the {\em sum of all increasing subsequences} is $2 \sqrt{2+\sqrt{5}}\, \sqrt{n}$.
One may easily optimize the formula above by solving for the critical point. It is 
$\kappa = \sqrt{2+2\sqrt{5}}/\sqrt{5}$ (which corresponds to $P(\kappa) = (\sqrt{5}-2)/4$). That does gives $\ln(\mathbf{E}[Z_{n,k_n}^2]) \sim 2 \sqrt{2+\sqrt{5}}\, \sqrt{n}$.

\section{An exactly solvable model: the number of small sums}
\label{sec:solvable}

Suppose that we have random variables $X_1,X_2,\dots$ which are IID, exponential-1 random variables:
\begin{equation}
	\forall n \in \{1,2,\dots\}\, ,\ \forall t_1,\dots,t_n \geq 0\, ,\qquad 
	\mathbf{P}(\{X_1\geq t_1\, ,\ \dots\, ,\ X_n\geq t_n\})\, =\, e^{-(t_1+\dots+t_n)}\, .
\end{equation}
Given $n \in \{1,2,\dots\}$ and $k \in [n]$ and given $t \geq 0$, let us define the random variable $\mathcal{N}_{n,k}(t)$ to be
\begin{equation}
	\mathcal{N}_{n,k}(t)\, =\, |\{ (i_1,\dots,i_k) \in [n]^k\, :\, i_1<\dots<i_k\ \text{ and }\ X_{i_1}+\dots+X_{i_k} \leq t\}|\, .
\end{equation}
We are going to consider $k$ depending on $n$ in such a way that $k/n^{1/2} \to \kappa \in (0,\infty)$.
Then we claim that the condition that $X_{i_1}+\dots+X_{i_k} \leq t$ has a resemblance to the analogous condition for random permutations
to have an increasing sequence.
This is most directly seen by noting that for any fixed $i_1<\dots<i_k$, the random variable $X_{i_1}+\dots+X_{i_k}$ is a Gamma-$(k,1)$ random variable.
Therefore, we have
\begin{equation}
	\mathbf{P}(\{X_{i_1}+\dots+X_{i_k} \leq t\})\, =\, \frac{1}{\Gamma(k)}\, \int_0^t x^{k-1} e^{-x}\, dx\, =\, \frac{\Gamma(k,t)}{\Gamma(k)}\, .
\end{equation}
Hence, we have
\begin{equation}
	\mathbf{E}[\mathcal{N}_{n,k}(t)]\, =\, \binom{n}{k} \, \cdot \frac{1}{\Gamma(k)}\, \int_0^t x^{k-1} e^{-x}\, dx,
\end{equation}
which implies that for a sequence $k_n$ such that $k_n \sim \kappa n^{1/2}$ as $n \to \infty$, we have
\begin{equation}
\lim_{n \to \infty} n^{-1/2} 
\ln\left(\mathbf{E}\left[\mathcal{N}_{n,k}(t)\right]\right)\, =\, \kappa\cdot \left(2 - 2 \ln(\kappa) + \ln(t)\right)\, .
\end{equation}
This may be compared to the first moment calculation for the generalized Ulam problem (\ref{eq:FirstMomentAsymp}).
But there the comparisons stop, because the correlations in the present problem are different.

\subsection{A replica symmetric ansatz}

Suppose that we have $m$ subsets $A_1,\dots,A_m \subset [n]$ each of cardinality $k$.
For each subset $R \subseteq [m]$ such that $R \neq \emptyset$, let us define the subset
\begin{equation}
	\mathcal{B}_{R}\, =\, \left\{i \in [n]\, :\, \left(\forall r \in R\, ,\ i \in A_r\right)\ \text{ and }\ \left(\forall r \in [m]\setminus R\, ,\ i \notin A_r\right)\right\}\, .
\end{equation}
Let us denote $\beta_R = |\mathcal{B}_{R}|$ to be the cardinality. And we observe that
\begin{equation}
	\forall r \in [m]\, ,\ \sum_{\substack{R \subseteq [m]\\ R \ni r}} \beta_R\, =\, |A_r|\, =\, k\, .
\end{equation}
Let us also define $\mathcal{T}_R = \sum_{i \in \mathcal{B}_R} X_i$.
So the condition that we wish to enforce is that 
\begin{equation}
	\left(\forall r \in [m]\, ,\ \sum_{i \in A_R} X_i \leq t\right)\qquad \Leftrightarrow\qquad
	\left(\forall r \in [m]\, ,\ \sum_{\substack{R \subseteq [m]\\ R \ni r}} \mathcal{T}_R \leq t\right)\, .
\end{equation}
Now the replica symmetric ansatz we make is that $\beta_R$ is only a function of $|R|$: $\forall  \ell \in [m]$, and for all $R \subseteq [m]$ such that $|R|=\ell$ we have $\beta_R=b_{\ell}$ for some fixed numbers $b_1,\dots,b_{\ell}$ satisfying
\begin{equation}
\sum_{\ell=1}^m \binom{m-1}{\ell-1} b_{\ell}\, =\, k\, .
\end{equation}
The reasons for this is that once $r \in [m]$ is selected, the number of subsets $R \subseteq [m]$ satisfying $|R|=\ell$ and $R \ni r$ is equal in number to $\binom{m-1}{\ell-1}$.
For symmetry reasons, we also assume that there are numbers $\tau_1,\dots,\tau_{m}$ such that 
\begin{equation}
\label{eq:constraint1}
\sum_{\ell=1}^m \binom{m-1}{\ell-1} \tau_{\ell}\, =\, t\, ,
\end{equation}
and we enforce (i.e., calculate the probability or more precisely the large deviation rate function for) 
\begin{equation}
	\forall R \in [m] \text{ (such that $R\neq \emptyset$), }\ \mathcal{T}_R \leq \tau_{|R|}\, .
\end{equation}
Then, since there are $\binom{m}{\ell}$ sets $R \subseteq [m]$ with $|R|=\ell$, we have, using Stirling's formula to evaluate the asymptotics of the multinomial coefficient,
\begin{equation}
\label{eq:Pasymp}
\begin{split}
\lim_{n \to \infty} n^{-1/2} \ln\left(
	\sum_{\substack{\{\mathcal{B}_R : R \subseteq [m]\, R\neq \emptyset\}\\ \forall \emptyset\neq R \subseteq [m]\, ,\ |\mathcal{B}_R|=b_{|R|}}} 
\mathbf{P}\left(\left\{\forall \emptyset\neq R \subseteq [m]\, ,\ \mathcal{T}_R \leq \tau_{|R|}\right\}\right)\right)\\
&\hspace{-4cm}=\, \sum_{\ell=1}^{m} \binom{m}{\ell} \kappa_{\ell} \left(2-2\ln(\kappa_{\ell})+\ln(\tau_{\ell})\right)\, ,
\end{split}
\end{equation}
where we assume $b_{\ell} \sim \kappa_{\ell} n^{1/2}$ for numbers $\kappa_{\ell}$ satisfying
\begin{equation}
\label{eq:constraint2}
\sum_{\ell=1}^m \binom{m-1}{\ell-1} \kappa_{\ell}\, =\, \kappa\, .
\end{equation}
Using Lagrange multipliers for the two constraints (\ref{eq:constraint1}) and (\ref{eq:constraint2}) it is easy to see that the critical points are such that
\begin{equation}
	\kappa_{\ell}\, =\, \frac{\sqrt{mt}}{\sqrt{f_m(z)}}\, \cdot \frac{z^\ell}{\ell}\qquad \text{ and }\qquad
	\tau_{\ell}\, =\, \frac{mt}{f_m(z)}\, \cdot \frac{z^{\ell}}{\ell^2}\, ,
\end{equation}
where $z$ solves
\begin{equation}
\label{eq:zImplicit}
	\frac{(1+z)^m-1}{\sqrt{f_m(z)}}\, =\, \frac{\kappa \sqrt{m}}{\sqrt{t}}\, ,
\end{equation}
for the function
$f_m(x)$ defined to be $\int_0^x \frac{(1+y)^m-1}{y}\, dy$ which may also be written as $\sum_{\ell=1}^{m} \binom{m}{\ell}\, \frac{x^\ell}{\ell}$.
Substituting this into (\ref{eq:Pasymp}), we obtain the formula
\begin{equation}
	\lim_{n \to \infty} n^{-1/2} \ln\left(\mathbf{E}\left[\mathcal{N}_{n,k}^{m}(t)\right]\right)\,  \text{\Huge ``}=\text{\Huge ''}\, 
\sqrt{\frac{mt}{f_m(z)}}\, \left[2 f_m(z) - \ln(z)\left((1+z)^m-1\right)\right]\, .
\end{equation}
Note that this formula is only an expectation predicated on the belief that the optimal choice will be symmetric under interchanges of the sets $A_1,\dots,A_m$.
In other words, we have made the replica symmetric ansatz. 

\subsection{The replica-to-zero trick}
If we extend the definition of $f_m(x)$ to more general parameters $m$, say called $\epsilon$, we obtain
$$
f_{\epsilon}(x)\, =\, \int_0^x \frac{e^{\epsilon \ln(1+y)} - 1}{y}\, dy\, \sim\, \epsilon \int_0^x \frac{\ln(1+y)}{y}\, dy\, =\, -\epsilon \operatorname{Li}_2(-x)\, ,
$$
as $\epsilon \to 0^+$,
using Spence's dilogarithm function $\operatorname{Li}_2(x)$ which is $ -\int_0^{x} \frac{\ln(1-x)}{x}\, dx$ or $\sum_{\ell=1}^{\infty} \frac{x^\ell}{\ell^2}$.
Then equation (\ref{eq:zImplicit}) becomes
$$
\frac{\ln(1+z)}{\sqrt{-\operatorname{Li}_2(-z)}}\, =\, \frac{\kappa}{\sqrt{t}}\, .
$$
We have the replica symmetric prediction from the replica-to-0 trick
$$
	\lim_{n \to \infty} n^{-1/2}\, \mathbf{E}\left[\ln\left(\mathcal{N}_{n,k}(t)\right)\right]\, \text{\Huge ``}=\text{\Huge ''}\,
	\sqrt{\frac{t}{-\operatorname{Li}_2(-z)}}\, \left[-2\operatorname{Li}_2(-z)-\ln(z)\ln(1+z)\right]\, .
$$
Here we use van Hemmen and Palmer's convention from \cite{vanHemmenPalmer}.

\subsection{Exact solution to leading order: validity of the replica symmetric ansatz}

The formula above can be seen to be incorrect, at least if $\kappa/\sqrt{t}$ is sufficiently large.
The distribution of the order statistics of $X_1,\dots,X_n$, near the minimum, is of a Poisson point process with spacing parameter $1/n$.
In other words,
$$
\lim_{n \to \infty} \mathbf{P}\left(\left\{|\{i \in [n]\, :\, X_i \leq t/n\}|=k\right\}\right)\, =\, \lim_{n \to \infty} \binom{n}{k}\, (1-e^{-t/n})^k e^{-t(n-k)/n}\, =\, e^{-t}\, \frac{t^k}{k!}\, ,
$$
for each $t>0$ and $k \in \{0,1,\dots\}$. By a similar calculation,
$$
\lim_{n \to \infty} n^{-1/2} \ln\left(\mathbf{P}\left(\left\{|\{i \in [n]\, :\, X_i \leq t/\sqrt{n}\}|=k\right\}\right)\right)\Bigg|_{k\sim \kappa \sqrt{n}}\, 
=\, \kappa\, \ln\left(\frac{t}{\kappa}\right) - t + \kappa\, ,
$$
which is maximized, as a function of $t$, at $t=\kappa$ where the limit is $0$.
This is what one expects for a Poisson point process, and it shows that the points are fairly tightly allocated to the lattice $\{1/n,2/n,\dots\}$, which is their collection
of their separate expectations.
But then asking for the number of ways of writing a particular order-1 number $t$ as a sum of such parts should be, to leading order, the same as the 
question of the number of partitions of $nt$ into distinct integers with approximately $\kappa \sqrt{n}$ parts.

This is a famous problem, apparently first considered by Szekeres using the circle method of Hardy and Ramanujan.
But much progress has been made to simplify the exposition, and we highly recommend the article of Dan Romik \cite{Romik}.
In particular, there are not any such partitions if $\kappa^2/2 \geq t$. In other words, we cannot have $\kappa/\sqrt{t} \geq \sqrt{2}$.
In fact, $\lim_{z \to \infty} \ln(1+z) / \sqrt{-\operatorname{Li}_2(-z)}=\sqrt{2}$. See, for example, the first reflection principle in Section 2 of Don Zagier's excellent 
reference \cite{Zagier}.

The formula appears to be correct for $z<1$ by using the the generating function for partitions
into distinct parts
$$
\prod_{n=1}^{\infty} (1+q^nz)\, =\, 1 + \sum_{k=1}^{\infty} \sum_{n=k(k+1)/2}^{\infty} \rho(n,k) q^n z^k\, ,
$$
where $\rho(n,k)$ is the number of partitions of the number $n$ with $k$ parts, all distinct. Then one may apply the Hayman saddle point method to the formula
$$
\rho(n,k)\, =\, \oint_{\mathcal{C}(0;a)} \oint_{\mathcal{C}(0;b)} \frac{e^{\sum_{n=1}^{\infty} \ln(1+q^nz)}}{q^n z^k}\, \cdot \frac{dq}{2\pi i q}\, \cdot \frac{dz}{2\pi i z}\, ,
$$
for $a=|q|<1,b=|z|<1$. One should take $n = t \sqrt{N}$ and $k = \kappa \sqrt{N}$ and choose the amplitude of the circle for $q$ such that when $\theta=0$ on the positive real axis
we have $q=e^{-\alpha/\sqrt{N}}$. Approximating the infinite sums by their Riemann integrals, and choosing $\alpha = \sqrt{-\operatorname{Li}_2(-z)/t}$
seems to give the correct formula for $|z|<1$.
For $|z|\geq 1$ some of the terms in the summation must be rewritten, because the exponential of the logarithm passes through a branch point.
But we do know that Romik's formula is already rigorously proved, and it is analytic with no phase transition. So the analytic extension
of the replica-to-zero formula also work here: one may simply take the analytic extension of the replica-to-$0$ formula from $|z|<1$ to $|z|\geq 1$.

Some technicaly details need to be verified to check the Hayman method for this problem, 
but we believe that it is correct.
It is not exactly the same presentation of the formula as in Romik's paper \cite{Romik}.
But we believe that the equivalence of the two formulas is due to a non-trivial identity for the dilogarithm function. 
Authors from Szekeres to Romik seem to use the identity $\rho(n,k) = p(n-\binom{k}{2},k)$, instead.
Here $p(n,k)$ is the usual
partition number: the number of partitions of $n$ into $k$ parts, not necessarily distinct.
(See for example, Theorem 3(b) on page 113 of Riordan \cite{Riordan}.)

We would like to note here that connections between mathematical
physics and topics such as asymptotic combinatorics or analytic number theory are one of the great dividends of mathematical physics.
One favorite example is that Vladimir Korepin has shown connections between the Riemann zeta function and the XXZ quantum spin system. See for example \cite{KorepinEtAl}.

Some authors study partitions as a stand-in for permutations or study permutations as a stand-in for primes, so that there is a close
connection between asymptotic combinatorics and analytic number theory. (For example, see the introduction to \cite{ABT}.)

\section{The third moment and positive RW's in $d=3$}

\label{sec:higher}

By Holder's inequality, we know
\begin{equation}
\label{eq:Holder}
\frac{1}{3\kappa\, \sqrt{n}}\, \ln\left(\mathbf{E}\left[Z_{n,k}^3\right]\right)\, \geq\, 
\frac{1}{2\kappa\, \sqrt{n}}\, \ln\left(\mathbf{E}\left[Z_{n,k}^2\right]\right)\, ,
\end{equation}
and we have precise asymptotics for the right-hand-side.
But we anticipate that the inequality is strict. One might hope that by correctly identifying the moments for a few powers that one could then guess the general
pattern and deduce the replica-to-0 limit.

For the moments higher than 2, there is an analogue of the combinatorial formula of Pinsky, but it does not necessarily capture the full picture.

\subsection{The all-or-nothing sector}
Here we consider $r$-tuples of subsequences such that if any two subsequences has a point in common, then that point is shared by all $r$ subsequences.
This simplifies calculations, and does provide a variational ansatz, although not necessarily the optimal trial state.
Then we generalize Pinsky's combinatorial array as follows:
\begin{equation}
	\widetilde{\mathcal{A}}_r(k_1,\dots,k_r,j)\, =\, \sum_{\substack{\alpha_1(0),\dots,\alpha_1(j) \in \{0,1,\dots\}\\ \alpha_1(0)+\dots+\alpha_1(j)=k_1}}
\quad \dots \quad 
 \sum_{\substack{\alpha_r(0),\dots,\alpha_r(j) \in \{0,1,\dots\}\\ \alpha_r(0)+\dots+\alpha_r(j)=k_j}} 
\prod_{s=0}^j \binom{\alpha_1(s)+\dots+\alpha_r(s)}{\alpha_1(s),\dots,\alpha_r(s)}^2\, . 
\end{equation} 
The second one is $\widetilde{\mathcal{A}}_2(k,\ell,j) = \mathcal{A}(k,\ell,j)$, the
original Pinsky formula.
Then we have the rigorous lower bound:
\begin{equation}
\label{ineq:RepSym}
	\mathbf{E}[Z_{n,k}^r]\, \geq\, \sum_{j=0}^{k} \mathbf{E}[Z_{n,rk-(r-1)j}] \widetilde{\mathcal{A}}_r(k-j,\dots,k-j,j)\, .
\end{equation}
This only holds for $r>1$. Note that for $r=2$ it is an identity because that is the original Pinsky formula.
For $r=1$ the only logical choices for $j$ are either: $j=k$ if you think of the intersections of the subsequence with itself; or else $j=0$ if you think of the intersections with an empty
second subsequence. So the sum would need to be replaced by just one of those two terms.
Either answer gives the correct formula for $r=1$ as an identity.


In order to analyze the asympotics, one needs the generating function. But first, 
let us define a generating function for the square of the multinomial coefficients $M^{(2)}_r(z_1,\dots,z_r)$:
\begin{equation}
	M^{(2)}_r(z_1,\dots,z_r)\, =\, \sum_{\alpha(1),\dots,\alpha(r) = 0}^{\infty} \binom{\alpha(1)+\dots+\alpha(r)}{\alpha(1),\dots,\alpha(r)}^2 
	z_1^{\alpha(1)} \cdots z_r^{\alpha(r)}\, .
\end{equation}
Then  the generating function for the original combinatorial quantity $\widetilde{\mathcal{A}}_r(k_1,\dots,k_r,j)$ satisfies
\begin{equation}
	\sum_{k_1,\dots,k_r,j = 0}^{\infty} \widetilde{\mathcal{A}}_r(k_1,\dots,k_r,j) x_1^{k_1} \cdots x_r^{k_r} w^j\,
	=\, \frac{1}{M^{(2)}_r(x_1,\dots,x_r)^{-1} - w}\, .
\end{equation}
From previous results  $M^{(2)}_2(x,y) = \left(1 - 2 x - 2 y + x^2 +y^2 - 2xy\right)^{-1/2}$.
Then, setting $y$ to be equal to $0$, we recover $M^{(1)}(x) = M^{(2)}(x,0) = 1/(1-x)$.
For later reference, let us rewrite in a slightly different way
\begin{equation}
\label{eq:M2Form}
	M^{(2)}_2(z_1,z_2)\, =\, \frac{1}{\Big(\big(1 - (z_1+z_2)\big)^2  - 4z_1z_2\Big)^{1/2}}\, .
\end{equation}

\subsection{The diagonal method for generating functions}

If we changed the square of the multinomial coefficient to the first power then the generating function would be trivial
\begin{equation}
\label{eq:M1formula}
	M^{(1)}_r(z_1,\dots,z_r)\, =\, \sum_{\alpha(1),\dots,\alpha(r) = 0}^{\infty} \binom{\alpha(1)+\dots+\alpha(r)}{\alpha(1),\dots,\alpha(r)} 
	z_1^{\alpha(1)} \cdots z_r^{\alpha(r)}\, =\, \frac{1}{1-(z_1+\dots+z_r)}\, .
\end{equation}
This generating function is valid as long as $|z_1|+\dots+|z_r|<1$, for example.
Then, as long as all the variables are in this domain, we may use Cauchy's integral formula to determine
\begin{equation}
\begin{split}
\oint_{\mathcal{C}(0;1)}\ \dots\ \oint_{\mathcal{C}(0;1)} M^{(1)}_r(x_1 \omega_1,\dots,x_r \omega_r)
\cdot M^{(1)}_r\left(\frac{y_1}{\omega_1},\dots,\frac{y_r}{\omega_r}\right)\, \frac{d\omega_1}{2\pi i \omega_1}\, \cdots\, 
\frac{d\omega_r}{2\pi i \omega_r}\\
&\hspace{-2cm} =\, M^{(2)}_r(x_1 y_1,\dots,x_r y_r)\, .   
\end{split}
\end{equation}
This is the ``diagonal method.''
Substituting the formula from (\ref{eq:M1formula}), we have
\begin{equation}
\oint_{\mathcal{C}(0;1)}\ \dots\ \oint_{\mathcal{C}(0;1)} \frac{1}{1-\sum_{s=1}^r x_s\omega_s} \cdot 
\frac{1}{1-\sum_{s=1}^{r} \frac{y_s}{\omega_s}}\, \cdot \frac{d\omega_1}{2\pi i \omega_1}\, \cdots\, 
\frac{d\omega_r}{2\pi i \omega_r}\,
=\, M^{(2)}_r(x_1 y_1,\dots,x_r y_r)\, .   
\end{equation}
This implies the formula
\begin{equation}
\label{eq:MGen}
\begin{split}
\oint_{\mathcal{C}(0;1)} \frac{1}{(1-x_r\omega_r)(\omega_r-y_r)}\, 
M^{(2)}_{r-1}\left(\frac{z_1 \omega_r}{(1-x_r\omega_r)(\omega_r-y_r)}\, ,\ \dots\, ,\ 
\frac{z_{r-1} \omega_r}{(1-x_r\omega_r)(\omega_r-y_r)}\right)\, \frac{d\omega_r}{2\pi i}\\
&\hspace{-3.5cm}
=\, M_r^{(2)}(z_1,\dots,z_{r-1},x_r y_r)\, ,
\end{split}
\end{equation}
which is our main recurrence formula.

\subsection{The generating function for $r=3$}

Using equation (\ref{eq:MGen}) and (\ref{eq:M2Form}), we have
\begin{equation}
M_3^{(2)}(z_1,z_2,x_3y_3)\,
=\, 
\oint_{\mathcal{C}(0;1)} \frac{1}{\left(\big((1-x_3\omega)(\omega-y_3)-  (z_1 + z_2) \omega\big)^2 - 4z_1 z_2\omega^2 \right)^{1/2}}\, \cdot
\frac{d\omega}{2\pi i}\, .
\end{equation}
For simplicity, let us introduce $x_1$ and $x_2$ and set $z_1=x_1^2$ and $z_2=x_2^2$. Let us also set $y_3=x_3$. So 
\begin{equation}
M_3^{(2)}(x_1^2,x_2^2,x_3^2)\,
=\, 
\oint_{\mathcal{C}(0;1)} \frac{1}{\left(\big((1-x_3\omega)(\omega-x_3)-  (x_1^2 + x_2^2) \omega\big)^2 - (2x_1x_2\omega)^2 \right)^{1/2}}\, \cdot
\frac{d\omega}{2\pi i}\, .
\end{equation}
Then the quantity inside the square-root is a difference of squares.
Defining the quartic in $\omega$ as 
\begin{equation}
	\mathcal{Q}(x_1,x_2,x_3;\omega)\, =\, \big((1-x_3\omega)(\omega-x_3)-  (x_1^2 + x_2^2) \omega\big)^2 - (2x_1x_2\omega)^2\, ,
\end{equation}
since it is a difference of squares, it may be factorized as a product of two quadratics
\begin{equation}
	\mathcal{Q}(x_1,x_2,x_3;\omega)\, =\, q_+(x_1,x_2,x_3;\omega)
\cdot  q_-(x_1,x_2,x_3;\omega)\, ,
\end{equation}
where
\begin{equation}
\begin{split}
	q_{\sigma}(x_1,x_2,x_3;\omega)\, 
&=\, (1-x_3\omega)(\omega-x_3)-  (x_1^2 + x_2^2) \omega + 2 \sigma x_1x_2 \omega\\
&=\, (1-x_3\omega)(\omega-x_3)-  (x_1 - \sigma x_2)^2 \omega\, .
\end{split}
\end{equation}
Then each quadratic equation may be factorized as 
\begin{equation}
	q_{\sigma}(x_1,x_2,x_3;\omega)\, 
=\, - x_3 \cdot \Big(\omega-\Omega_{\sigma,+}(x_1,x_2,x_3)\Big)
\cdot \Big(\omega-\Omega_{\sigma,-}(x_1,x_2,x_3)\Big)\, ,
\end{equation}
for 
\begin{equation}
\begin{split}
\Omega_{\sigma,\tau}(x_1,x_2,x_3)\,
&=\, \frac{\left(\sqrt{\big(1+x_3\big)^2-\big(x_1-\sigma x_2\big)^2} 
+ \tau \sqrt{\big(1-x_3\big)^2-\big(x_1-\sigma x_2\big)^2}\right)^2}{4x_3}\, .
\end{split}
\end{equation}
From the nature of the quadratics (having the same constant and leading coefficient) the roots have the property
$\Omega_{\sigma,+}(x_1,x_2,x_3) \Omega_{\sigma,-}(x_1,x_2,x_3) = 1$.
For $x_1,x_2,x_3$ small but positive, $\Omega_{\sigma,+}(x_1,x_2,x_3)$ will be positive and larger than 1,
while $\Omega_{\sigma,-}(x_1,x_2,x_3)$ will be positive and smaller than 1.
More precisely, for $x_1,x_2,x_3$ all small and positive, and of the same order, we have
$\Omega_{\sigma,-}(x_1,x_2,x_3) -x_3 \sim x_3 (x_1-\sigma x_2)^2$
and $\Omega_{\sigma,+}(x_1,x_2,x_3) -x_3^{-1} \sim -x_3^{-1} (x_1-\sigma x_2)^2$.
From this, we see that in this regime, we have
\begin{equation}
	0<\Omega_{+,-}(x_1,x_2,x_3) <\Omega_{-,-}(x_1,x_2,x_3) < 1 < \Omega_{-,+}(x_1,x_2,x_3) < \Omega_{+,+}(x_1,x_2,x_3)\, .
\end{equation}
Hence the $\omega$ contour integral about $\mathcal{C}(0;1)$ for the integrand $1/\sqrt{\mathcal{Q}(x_1,x_2,x_3;\omega)}$
can be deformed to an elliptic integral
\begin{equation}
M_3^{(2)}(x_1^2,x_2^2,x_3^2)\,
=\, \frac{1}{\pi}\, \int_{\Omega_{+,-}(x_1,x_2,x_3)}^{\Omega_{-,-}(x_1,x_2,x_3)} \frac{dt}{\sqrt{\widetilde{Q}(x_1,x_2,x_3;t)}}\, ,
\end{equation}
where the new quartic is positive on the interval of integration:
\begin{equation}
\begin{split}
\widetilde{Q}(x_1,x_2,x_3;t)\, &=\,
x_3^2 \cdot
\big(t-\Omega_{-,-}(x_1,x_2,x_3)\big)
\big(\Omega_{+,-}(x_1,x_2,x_3)-t\big)\\
&\qquad \qquad \cdot
\big(\Omega_{+,+}(x_1,x_2,x_3)-t\big)
\big(\Omega_{-,+}(x_1,x_2,x_3)-t\big)\, .
\end{split}
\end{equation}
This may be transformed into the standard form for a complete elliptic integral by applying a linear fractional transformation.
The Legendre form of the complete elliptic integral of the 1st kind is 
\begin{equation}
	K(k)\, =\, \int_0^1 \frac{dt}{\sqrt{(1-t^2)(1-k^2t^2)}}\, 
=\, \frac{1}{2}\, \int_{-1}^{1} \frac{dt}{\sqrt{(1-t^2)(1-k^2t^2)}}\, ,
\end{equation}
for $0\leq k<1$. (Note that this version of the elliptic integral is itself the generating function for the probability of return of a 2d SRW
to the origin after $2n$ steps [which is easiest to see using the diagonal transformation that is also the key to Pinsky's combinatorial 2nd moment formula]. See for example P\'olya \cite{Polya}. Or for a more accessible alternative see
McKean’s graduate probability textbook \cite{McKean} in which he discusses
P\'olya’s result in Section 4.3.)
Given this, we may write
\begin{equation}
\label{eq:MFormula}
	M_3^{(2)}(x_1^2,x_2^2,x_3^2)\, =\, m_3^{(2)}(x_1,x_2,x_3) K(k(x_1,x_2,x_3))\, ,
\end{equation}
where we define
\begin{equation}
	\varkappa_{\tau}(x_1,x_2,x_3)\, =\, \prod_{\substack{\sigma_1,\sigma_2,\sigma_3 \in \{+1,-1\}\\ \sigma_1\sigma_2\sigma_3=\tau}} \sqrt{1+\sigma_1 x_1 +\sigma_2 x_2 + \sigma_3 x_3}\, ,\qquad
\text{ for $\tau \in \{+1,-1\}$,}
\end{equation}
and then we may write
\begin{equation}
\label{eq:kFormula}
	k(x_1,x_2,x_3)\, =\, \left(\frac{\sqrt{\varkappa_+(x_1,x_2,x_3)} - \sqrt{\varkappa_-(x_1,x_2,x_3)}}
{\sqrt{\varkappa_+(x_1,x_2,x_3)} + \sqrt{\varkappa_-(x_1,x_2,x_3)}}\right)^2\, ,
\end{equation}
and
\begin{equation}
	m_2^{(2)}(x_1,x_2,x_3)\, =\, \frac{8}{\pi \cdot \left(\sqrt{\varkappa_+(x_1,x_2,x_3)} + \sqrt{\varkappa_-(x_1,x_2,x_3)}\right)^2}\, .
\end{equation}
One positive aspect of this format is that this is clearly symmetric in $x_1,x_2,x_3$, as it should be.

\subsubsection{Specialization to $x_3=0$ recovers $M_2^{(2)}(x,y)$}

Taking the $x_3 \to 0^+$ limit, we must recover $M_3^{(2)}(x_1^2,x_2^2,0) = M_2^{(2)}(x_1^2,x_2^2)$.
If we set $x_3=0$ then symmetry implies $\varkappa_+(x_1,x_2,0)=\varkappa_-(x_1,x_2,0)$.
More precisely,
\begin{equation}
	\varkappa_{+}(x_1,x_2,0)\, =\, \varkappa_{-}(x_1,x_2,0)\, =\, 
	\prod_{\sigma_1,\sigma_2 \in \{+1,-1\}} \sqrt{1+\sigma_1 x_1 +\sigma_2 x_2}\, ,
\end{equation}
because the sign $\sigma_3$ for $x_3$ may be either $+1$ or $-1$, as desired, without affecting the argument: $+0=-0=0$.
Moreover, this quantity is $\sqrt{1+x_1+x_2}\, \sqrt{1-x_1-x_2}\, \sqrt{1+x_1-x_2}\, \sqrt{1-x_1+x_2}$.
This may be written more simply as a product of differences of squares
\begin{equation}
	\varkappa_{\tau}(x_1,x_2,0)\, =\, \sqrt{1-(x_1+x_2)^2}\, \sqrt{1-(x_1-x_2)^2}\, .
\end{equation}
Then taking the product, by bringing both terms to the same square-root, we have
\begin{equation}
	\varkappa_{\tau}(x_1,x_2,0)\, =\, \sqrt{1-2(x_1^2+x_2^2) - (x_1^2-x_2^2)^2}\, .
\end{equation}
Then writing $z_1=x_1^2$ and $z_2=x_2^2$ we have 
\begin{equation}
	\varkappa_{\tau}(x_1,x_2,0)\,  =\, \sqrt{1-2(z_1+z_2) - (z_1-z_2)^2}\, .
\end{equation}
So this gives
\begin{equation}
\begin{split}
	m_2^{(2)}(x_1,x_2,x_3)\, &=\, \frac{8}{\pi \cdot \left(\sqrt{\varkappa_+(x_1,x_2,0)} + \sqrt{\varkappa_-(x_1,x_2,0)}\right)^2}\\
	&=\, \frac{2}{\pi\, \sqrt{1-2(x_1^2+x_2^2) - (x_1^2-x_2^2)^2}}\, .
\end{split}
\end{equation}
Since $\varkappa_+(x_1,x_2,0)=\varkappa_-(x_1,x_2,0)$, the formula in (\ref{eq:kFormula}) becomes $k(x_1,x_2,0)=0$.
But the elliptic integral at $0$ is just 
\begin{equation}
	K(0)\, =\, \int_0^1 \frac{dt}{\sqrt{1-t^2}}\, =\, \frac{\pi}{2}\, .
\end{equation}
So we do have from (\ref{eq:MFormula}) that
\begin{equation}
	M_3^{(2)}(z_1,z_2,0)\, =\, \frac{1}{\sqrt{1-2(z_1+z_2) - (z_1-z_2)^2}}\, ,
\end{equation}
as in (\ref{eq:M2Form}).

\subsection{Some observations}

This generating function has value on its own. One may think of the 2nd moment generating function that we calculated before as characterizing the number
of intersections of two up-right paths. By rotating 45-degrees and projecting to the $x$-axis, this can be seen as the rencontres of two SRW's on $\Z$.
Those are the times that both SRW's end up at the same point at one time.
For the 3rd moment, the original perspective is of $(+\boldsymbol{i},+\boldsymbol{j},+\boldsymbol{k})$ paths in $\Z^3$.
This can be thought of as positive RW's in $d=3$.
If we project to the plane $x+y+z=0$ along the $(1,1,1)$ direction, this can be seen as a particular type of RW along the positive directions of the triangular lattice.

Another observation is that the all-or-nothing sector does not lead to a better inequality
than the one we started with (\ref{eq:Holder}). This is discovered {\em a posteriori} after doing much calculation.
Therefore, we will skip those unnecessary calculations.

\section{Proofs of the combinatorial identities}
\label{sec:Pictures}

The first proof is a direct adaptation of Pinsky's proof of his Proposition 2 in \cite{Pinsky1},
just as the statement of the proposition is a direct adaptation of his discovery.

\begin{proofof}{\bf Proof of Propsosition \ref{prop:PinskyGeneralization}:}
Let $\mathcal{C}_n(k,\ell,j)$ be the set of all pairs $(\boldsymbol{a},\boldsymbol{b})$ where 
$\boldsymbol{a} = (a_1,\dots,a_k) \in [n]^k$ and $\boldsymbol{b} \in (b_1,\dots,b_{\ell}) \in [n]^{\ell}$ are tuples satisfying
$$
a_1<\dots<a_k\ \text{ and }\ b_1<\dots<b_{\ell}\, ,\ \text{ and also }\
|\{a_1,\dots,a_k\} \cap \{b_1,\dots,b_{\ell}\}|\, =\, j\, .
$$
Let us be more precise.
Let $\Lambda_j(k-j,\ell-j)$ denote the set of all pairs $(\boldsymbol{\alpha},\boldsymbol{\beta})$
where $\boldsymbol{\alpha} = (\alpha_0,\dots,\alpha_j)$
and $\boldsymbol{\beta} = (\beta_0,\dots,\beta_j)$ are tuples in
$\{0,1,\dots\}^{j+1}$ satisfying
$$
\alpha_0+\dots+\alpha_j\, =\, k-j\ \text{ and }\ 
\beta_0+\dots+\beta_j\, =\, \ell-j\, .
$$
Then, given $(\boldsymbol{\alpha},\boldsymbol{\beta}) \in \Lambda_j(k-j,\ell-j)$, let 
us define $\mathcal{C}_n(k,\ell,j;\boldsymbol{\alpha},\boldsymbol{\beta})$
to be the subset consisting of all pairs $(\boldsymbol{a},\boldsymbol{b}) \in \mathcal{C}_n(k,\ell,j)$
satisfying the following further condition. Given $(\boldsymbol{\alpha},\boldsymbol{\beta}) \in \Lambda_j(k-j,\ell-j)$,
let us define $M(1),\dots,M(j)$ and $N(1),\dots,N(j)$ such that we first set $M(0)=N(0)=0$ and then define
for each $r \in \{1,\dots,j\}$,
$$
M(r)\, =\, M(r-1) + \alpha_r+1\ \text{ and }\ N(r)\, =\, N(r-1)+\beta_r+1\, .
$$
Then the condition for $(\boldsymbol{a},\boldsymbol{b}) \in \mathcal{C}_n(k,\ell,j)$
to also be in $\mathcal{C}_n(k,\ell,j;\boldsymbol{\alpha},\boldsymbol{\beta})$ is that
$$
a_{M(r)}\, =\, b_{M(r)}\ \text{ for all $r \in \{1,\dots,j\}$.}
$$
Of course, since $|\{a_1,\dots,a_k\}\cap \{b_1,\dots,b_k\}|=j$, that means that these are exactly the points of intersection.
Now we will prove the following
$$
	\sum_{(\boldsymbol{a},\boldsymbol{b}) \in \mathcal{C}_n(k,\ell,j;\boldsymbol{\alpha},\boldsymbol{\beta})}
	\mathbf{P}(\pi_{a_1}<\dots<\pi_{a_k}\ \text{ and }\ \pi_{b_1}<\dots<\pi_{b_{\ell}})\,
=\, \binom{n}{k+\ell-j} \frac{1}{(k+\ell-j)!}\, \prod_{r=0}^{j} \binom{\alpha_r+\beta_r}{\alpha_r,\beta_r}^2\, .
$$
This formula would establish the proof we are seeking.
We show this in two steps, corresponding to the two powers of the $\alpha$-$\beta$ binomials.

Given $(\boldsymbol{a},\boldsymbol{b})$, let $\boldsymbol{c}=(c_1,\dots,c_{k+\ell-j}) \in [n]^{k+\ell-j}$ be the tuple 
such that
$$
c_1<\dots<c_{k+\ell-j}\ \text{ and }\ 
\{c_1,\dots,c_{k+\ell-j}\}\, =\, 
\{a_1,\dots,a_k\} \cup \{b_1,\dots,b_{\ell}\}\, .
$$
Likewise, let $\boldsymbol{d}=(d_1,\dots,d_j) \in [n]^{k+\ell-j}$ be the tuple such that
$$
d_1<\dots<d_{j}\ \text{ and }\ 
\{d_1,\dots,d_{j}\}\, =\, 
\{a_1,\dots,a_k\} \cap \{b_1,\dots,b_{\ell}\}\, .
$$
Assuming $(\boldsymbol{\alpha},\boldsymbol{\beta})$ is the pair such that
$(\boldsymbol{a},\boldsymbol{b}) \in \mathcal{C}_n(k,\ell,j;\boldsymbol{\alpha},\boldsymbol{\beta})$,
we may note the following.
Firstly, we may define $P(1),\dots,P(j)$ such that first we set $P(0)=0$ and then define for each $r \in \{1,\dots,j\}$
$$
P(r)\, =\, P(r-1)+\alpha_r+\beta_r+1\, .
$$
Secondly, we then have
$$
	c_{M(r)}\, =\, d_r\ \text{ for all }\ r \in \{1,\dots,j\}\, .
$$
Now we will first prove that
\begin{equation}
\label{eq:ProbabEqual}
\mathbf{P}(\pi_{a_1}<\dots<\pi_{a_k}\ \text{ and }\ \pi_{b_1}<\dots<\pi_{b_{\ell}})\,
=\, \mathbf{P}(\pi_{c_1}<\dots<\pi_{c_{k+\ell-j}})\, \prod_{r=0}^{\ell} \binom{\alpha_r+\beta_r}{\alpha_r,\beta_r}\, .
\end{equation}

To see this, let $\pi^{(1)} \in S_n$ be any permutation such that 
$$
\pi^{(1)}_{c_1}<\dots<\pi^{(1)}_{c_{k+\ell-j}}\, .
$$
Then we may construct another permutation $\pi^{(2)} \in S_n$ such that
$$
\pi^{(2)}_{a_1}<\dots<\pi^{(2)}_{a_k}\ \text{ and }\ \pi^{(2)}_{b_1}<\dots<\pi^{(2)}_{b_{\ell}}
$$
in multiple ways, starting from $\pi^{(1)}$, as we now describe.
Let $\pi^{(2)}_i = \pi^{(1)}_i$ for every $i \in [n] \setminus \{c_1,\dots,c_{k+\ell-j}\}$.
Let $\pi^{(2)}_{d_r} = \pi^{(1)}_{d_r}$ for every $r \in \{1,\dots,j\}$.
But for each $r \in \{0,\dots,j\}$, let us define $\pi^{(2)}_i$ for $i \in \{c_{P(r-1)+1},\dots,c_{P(r)-1}\}$
as follows. Choose any cardinality $\alpha_j$ subset of 
$$
\{\pi^{(1)}_{c_{P(r-1)+1}},\dots,\pi^{(1)}_{c_{P(r)-1}}\}\, ,
$$
to be the numbers (keeping their relative order) corresponding to 
$$
\{\pi^{(2)}_{a_{M(r-1)+1}},\dots,\pi^{(2)}_{a_{M(r)-1}}\}\, .
$$
Note that the complementary set will have cardinality $\beta_j$ and will be the number (keeping their relative order)
corresponding to 
$$
\{\pi^{(2)}_{b_{N(r-1)+1}},\dots,\pi^{(2)}_{b_{N(r)-1}}\}\, .
$$
So we claim that we have established (\ref{eq:ProbabEqual}).
Note that since $\mathbf{P}(\pi_{c_1}<\dots<\pi_{c_{k+\ell-j}})=1/(k+\ell-j)!$, we have
$$
\mathbf{P}(\pi_{a_1}<\dots<\pi_{a_k}\ \text{ and }\ \pi_{b_1}<\dots<\pi_{b_{\ell}})\,
=\, \frac{1}{(k+\ell-j)!}\, \prod_{r=0}^{\ell} \binom{\alpha_r+\beta_r}{\alpha_r,\beta_r}\, .
$$

Now we turn attention to the proof of the second part, which is that
$$
|\mathcal{C}_n(k,\ell,j;\boldsymbol{\alpha},\boldsymbol{\beta})|\, =\, |\mathcal{C}_n(k+\ell-j)|\, \prod_{r=0}^{j}
\binom{\alpha_r+\beta_r}{\alpha_r,\beta_r}\, ,
$$
where $\mathcal{C}_n(k+\ell-j)$ refers to all tuples $\boldsymbol{c} = (c_1,\dots,c_{k+\ell-j}) \in [n]^{k+\ell-j}$ satisfying
$$
c_1\, <\, \dots\, <\, c_{k+\ell-j}\, .
$$
Of course, then the cardinality of $\mathcal{C}_n(k+\ell-j)$ is equal to the binomial coefficient $n$-choose-$(k+\ell-j)$.
Thus, with the first part, this second part will establish the desired formula, which in turn will establish the proof of the proposition.
Given $\boldsymbol{c}^{(1)} \in \mathcal{C}_n(k+\ell-1)$ and given $j$ and $(\boldsymbol{\alpha},\boldsymbol{\beta})$, we choose $(\boldsymbol{a}^{(2)},\boldsymbol{b}^{(2)})$ as follows.
Firstly, we let $a^{(2)}_{M(r)}=b^{(2)}_{N(r)}=c^{(1)}_{P(r)}$.
Then for each $r$ we choose a cardinality $\alpha_r$ subset of 
$$
\{c^{(1)}_{P(r-1)+1},\dots,c^{(1)}_{P(r)-1}\}
$$
to be (in their relative order)
$$
(a^{(2)}_{M(r-1)+1},\dots,a^{(2)}_{M(r)-1})\, .
$$
Note that the complementary subset is cardinality $\beta_r$ and (in their relative order) corresponds to 
$$
(b^{(2)}_{N(r-1)+1},\dots,b^{(2)}_{N(r)-1})\, .
$$
\end{proofof}

We note that conditioning permutations on blocks and considering the projections of the points in 
those blocks to their $x$ coordinates and $y$ coordinates separately is also important for the Mallows
model on permutations, as in \cite{StarrWalters}.
But Pinsky's work predates that.

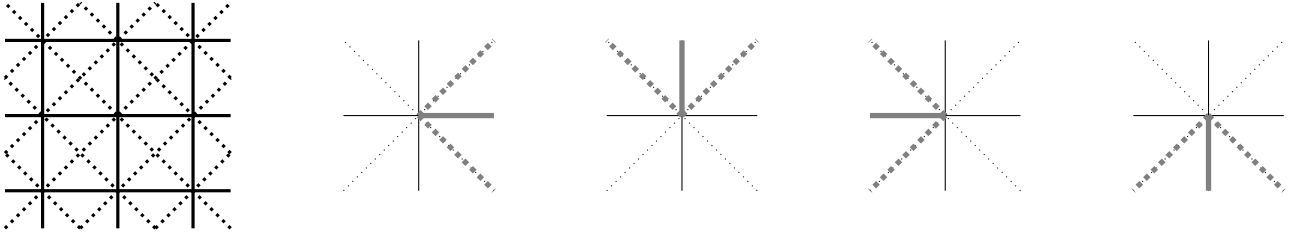
\begin{figure}
\begin{center}
\begin{tikzpicture}
\begin{scope}[very thick]
\foreach \x in {1,2,3}
{\draw[] (\x,0.5) -- (\x,3.5);}
\foreach \y in {1,2,3}
{\draw[] (0.5,\y) -- (3.5,\y);}
\draw[dotted] (0.5,2.5) -- (1.5,3.5);
\draw[dotted] (0.5,1.5) -- (2.5,3.5);
\draw[dotted] (0.5,0.5) -- (3.5,3.5);
\draw[dotted] (1.5,0.5) -- (3.5,2.5);
\draw[dotted] (2.5,0.5) -- (3.5,1.5);
\begin{scope}[xshift=4cm,xscale=-1]
\draw[dotted] (0.5,2.5) -- (1.5,3.5);
\draw[dotted] (0.5,1.5) -- (2.5,3.5);
\draw[dotted] (0.5,0.5) -- (3.5,3.5);
\draw[dotted] (1.5,0.5) -- (3.5,2.5);
\draw[dotted] (2.5,0.5) -- (3.5,1.5);
\end{scope}
\end{scope}
\begin{scope}[xshift = 5cm,yshift=1cm]
	\draw (0,1) -- (2,1);
	\draw (1,0) -- (1,2);
	\draw[dotted] (0,0) -- (2,2);
	\draw[dotted] (2,0) -- (0,2);
	\draw[line width=2pt,black!50!white,opacity=0.5] (1,1) -- (2,1);
	\draw[line width=2pt,black!50!white,opacity=0.5,dotted] (1,1) -- (2,2);
	\draw[line width=2pt,black!50!white,opacity=0.5,dotted] (1,1) -- (2,0);
\begin{scope}[xshift = 3.5cm]
	\draw (0,1) -- (2,1);
	\draw (1,0) -- (1,2);
	\draw[dotted] (0,0) -- (2,2);
	\draw[dotted] (2,0) -- (0,2);
	\draw[line width=2pt,black!50!white,opacity=0.5] (1,1) -- (1,2);
	\draw[line width=2pt,black!50!white,opacity=0.5,dotted] (1,1) -- (2,2);
	\draw[line width=2pt,black!50!white,opacity=0.5,dotted] (1,1) -- (0,2);
\begin{scope}[xshift = 3.5cm]
	\draw (0,1) -- (2,1);
	\draw (1,0) -- (1,2);
	\draw[dotted] (0,0) -- (2,2);
	\draw[dotted] (2,0) -- (0,2);
	\draw[line width=2pt,black!50!white,opacity=0.5] (1,1) -- (0,1);
	\draw[line width=2pt,black!50!white,opacity=0.5,dotted] (1,1) -- (0,2);
	\draw[line width=2pt,black!50!white,opacity=0.5,dotted] (1,1) -- (0,0);
\begin{scope}[xshift = 3.5cm]
	\draw (0,1) -- (2,1);
	\draw (1,0) -- (1,2);
	\draw[dotted] (0,0) -- (2,2);
	\draw[dotted] (2,0) -- (0,2);
	\draw[line width=2pt,black!50!white,opacity=0.5] (1,1) -- (1,0);
	\draw[line width=2pt,black!50!white,opacity=0.5,dotted] (1,1) -- (2,0);
	\draw[line width=2pt,black!50!white,opacity=0.5,dotted] (1,1) -- (0,0);
\end{scope}\end{scope}
\end{scope}
\end{scope}
\end{tikzpicture}
\caption{
The diagonal transformation that is part of Pinsky's combinatorial theorem. The choice of the next step being among $\{(0,1),(1,0),(-1,0),(0,-1)\}$ is equivalent
to a diagonal version of the encoding being in the set $\{(+1,+1),(+1,-1),(-1,-1),(-1,+1)\}$ as shown.
\label{fig:diag}
}
\end{center}
\end{figure}

\begin{proofof}{\bf Proof of Corollary \ref{cor:RW}:}
This proof follows Pinsky's argument.
Given two independent 1d random walks $X_t$ and $Y_t$
one may make a 2d random walk by taking $W_t = ((X_t+Y_t)/2,(X_t-Y_t)/2)$.
Geometrically, this corresponds to a certain diagonal transformation, that we  illustrate
in Figure \ref{fig:diag}.
Then two have $W_t$ on the $x$-axis for some $t$ requires that $X_t=Y_t$ for that $t$.

Let $t_0,t_1,\dots,t_j$ be a collection of times in non-decreasing order such that $X_{t_r}=Y_{t_r}$
for each $r \in \{0,1,\dots,j\}$. Assume $t_j=k+\ell$. Condition on the values $x_r = X_{t_r}$.
Let $\alpha_r$ be the number of up-steps and $\beta_r$ be the number of down-steps
of $X_t$, and hence also of $Y_t$, between times $t_{r-1}$ and $t_r$.
In other words, $\alpha_r-\beta_r=x_r-x_{r-1}$. (Let $t_{-1}=0$ and $x_{-1}=0$.)
Then the probability of these numbers of steps-up and steps-down at each inter-time block
for both $X_t$ and $Y_t$ is the product by independence $2^{-2k-2\ell}\prod_{r=0}^{j} \binom{\alpha_r+\beta_r}{\alpha_r,\beta_r}^2$.
Note that at time $t_j=k+\ell$ we have $x_j=\sum_{r=0}^{j} (\alpha_j-\beta_j)=k-\ell$.
Therefore, $W_{k+\ell} = (k-\ell,0)$.
\end{proofof}

\begin{proofof}{\bf Proof of Corollary \ref{cor:GF}:}
We begin with the proof of equation (\ref{eq:CombGF}).
By Newton's version of the binomial theorem,
\begin{equation*}
\begin{split}
\frac{1}{\left(1-2(x+y)+(x-y)^2\right)^{1/2}}\, 
&=\,  \frac{1}{1-(x+y)}\, \left(1- \frac{4xy}{\left(1-(x+y)\right)^2}\right)^{-1/2}\\
&=\, \sum_{n,k=0}^{\infty} \frac{(-1/2)_n}{n!} \cdot \frac{(-2n-1)_k}{k!} (-1)^{n+k} 4^n x^n y^n (x+y)^k\, ,
\end{split}
\end{equation*}
where the Pochhammer symbol, falling factorial, is 
\begin{equation*}
	(z)_n\, =\, z(z-1)(z-2)\cdots (z-n+1)\, =\, \prod_{k=0}^{n-1} (z-k)\, .
\end{equation*}
Then, by the regular binomial formula, we have
\begin{equation}
\frac{1}{\left(1-2(x+y)+(x-y)^2\right)^{1/2}}\, 
	=\, \sum_{n,k=0}^{\infty} \sum_{j=0}^{k} \frac{(-1/2)_n}{n!} \cdot \frac{(-2n-1)_k}{k!} \binom{k}{j} (-1)^{n+k} 4^n x^{n+j} y^{n+k-j}\, .
\end{equation}
Therefore, the equation reduces to the identity
\begin{equation*}
	\sum_{n=0}^{\infty} \frac{(-1/2)_n}{n!} \cdot \frac{(-2n-1)_k}{k!} \binom{k}{j} (-1)^{n+k} 4^n \Bigg|_{\substack{j=\ell-n\\k-j=m-n}}\,
	=\, \left(\binom{\ell+m}{\ell,m}\right)^2\, ,
\end{equation*}
for each $\ell,m \in \{0,1,\dots\}$.
By well-known properties of the Pochhammer symbols shifted by $1/2$, this is equivalent to checking
\begin{equation}
\label{eq:gamma2combinatorial}
	\sum_{n=0}^{\infty} \binom{\ell+m}{n,n,\ell-n,m-n}\, =\, \left(\binom{\ell+m}{\ell,m}\right)^2\, .
\end{equation}
But this follows from the previous lemma, using Pinsky's diagonal argument.
Suppose that after $\ell+m$ steps we know $W_{\ell+m} = (\ell-m,0)$.
Let $n$ be the number of left-steps, so that $n$ is also the number of right-steps.
Then there are $\ell+m-2n$ steps left. So the number of up-steps must be $\ell-n$ and the number
of down-steps must be $m-n$. Summing over $n$ we get all possible ways to have $W_{\ell+m}=(\ell-m,0)$.
But Pinsky's diagonal argument shows that this is the same as the number of ways of having
$X_{k+\ell}=Y_{k+\ell}=k-\ell$, for $X_t$ and $Y_t$ independent 1d random walks.
The number of these ways is the right-hand-side of (\ref{eq:gamma2combinatorial}).

The lemma follows from the elementary observation that 
\begin{equation}
	\sum_{k=0}^{\infty} \sum_{\ell=0}^{\infty} \sum_{j=0}^{\infty} \mathcal{A}(k,\ell,j) x^k y^{\ell} z^j\, =\, \sum_{j=0}^{\infty} \left(\sum_{\alpha=0}^{\infty} \sum_{\beta=0}^{\infty} \binom{\alpha+\beta}{\alpha,\beta}^2 x^{\alpha} y^{\beta}\right)^{j+1} z^j\, .
\end{equation}
\end{proofof}

\section{Discussion and outlook}

We have not proved every detail here.
The most challenging technical point is to check that the saddle point method applies. 
For the second moment and the associated combinatorial problem, we indicate the method in Appendices A through D below.
A good reference is the textbook of Pemantle and Wilson \cite{PemantleWilson}.

All other skipped proofs are algebraic or combinatorial and are much easier.
The proof of Lemma 2.5 is identical to the analysis of Subsections 5.1 and 5.2, specialized for $\mathcal{M}_2(x)=(1-4x)^{-1/2}$.
Theorem 2.4 is similar to the analysis of Subsection 5.1, properly generalized to the asymmetric case.
Then Theorem 2.6 is similar to the analysis of Subsection 5.2 in the asymmetric case.
(We include the long but elementary derivations in the arXiv preprint.)

For Section 4, the main issue is proper application of the diagonal method. Note that the simplest types
of inner-automorphisms for simplifying elliptic integrals are the linear-fractional transformations.
This is how to transform equation (37) to equation (40).
Since $\widetilde{Q}$ has two pairs of roots, where each pair consists of two roots whose product is 1,
the first good transformation is 
$$
\mathcal{L}(z)\, =\, \frac{z-1}{z+1}\, ,\qquad 
\mathcal{L}^{-1}(z)\, =\, \frac{z+1}{z-1}\, ,
$$
because $\mathcal{L}(1/z)=-\mathcal{L}(z)$. So it intertwines the multiplicative and additive inverses.
Then the second transformation uses the cross-ratio. More precisely, we want to take the new roots
$(-\rho_+,-\rho_-,\rho_-,\rho_+)$
to $(-k,k,1,-1)$ in that order, with a linear fractional transformation.
This works if and only if the cross ratios are equal $[-\rho_+,-\rho_-;\rho_-,\rho_+] = [-k,k;1,-1]$.
This leads to a formula for $k$, which is the one given in equation (42).
The appropriate linear fractional transformation is 
$$
\Lambda(z)\, =\, \sqrt{\rho_+ \rho_-}\, \cdot \frac{z+k^{1/2}}{-z+k^{1/2}}\, .
$$
Details beyond these are in the arXiv preprint, because they are long but elementary.

Now we mention some comments and open questions.

\begin{itemize}
\item In the diagonal method of Furstenberg and of Hautas and Klarner, it is known that rational functions
are mapped to algebraic functions after 1 step of the diagonal method.
See for example, the textbook of Pemantle and Wilson \cite{PemantleWilson}.
In the example of Section 4, the algebraic function $\mathcal{M}_2(x,y) = (1-2(x+y)+(x-y)^2)^{-1/2}$
is mapped to a function $\mathcal{M}_3(x,y,z)$ involving algebraic terms and a complete elliptic integral.
Is this typical? Also, what lies beyond, at the next step?
\item We wish to return to the consideration of the replica symmetric ansatz and replica-to-0 trick for the generalized Ulam problem.
This seems challenging, but worthwhile.
\item We are actually most interested in the step beyond the one we are stopped at.
Imagine any combinatorial stochastic optimization problem, as in the monograph of Steele \cite{Steele}.
Suppose that we can calculate all the moments for the number of elements whose cost is less than some number
$-k$, for each $k$.
How does one go from that to the probability that there exists at least one element whose 
cost is less than $-k$? The natural guess is the Bonferroni inequalities.
(See, for example, Chapter 4 of Charalambides \cite{Charalambides}.)
Looking at a generating function, it seems possible that the contour integral which arises is not 
Hayman-admissible. We are interested in that problem for a model where it applies, such as the solvable model we introduced.
Van Hemmen and Palmer proposed that replica symmetry breaking may be related to indeterminacy of the moment problem.
There, the ambiguity is characterized by a probability measure on the circle.
We are wondering whether the replica symmetric ansatz may fail when large-deviation saddle point calculations
associated to Bonferroni's inequalities require the full circle method of 
Hardy and Ramanujan.
\end{itemize}

Obviously a question that we are interested in is whether the replica method has anything interesting
to tell us about Ulam's problem, or vice-versa.
In general, the replica method is mysterious.
If one can apply it to more exactly solvable models, it might help to resolve some of the unknowns about why it works when it does.

\section*{Acknowledgments}

S.S.~is grateful to Pawel Hitczenko for a helpful discussion at RS\&A2023.
In particular, thanks for the reference to \cite{FlajoletSedgewick}.
The work of S.S.~was partly supported by 
a Simons collaboration grant.

\appendix 

%
%
%
%
%
%
%
%
%
%
%
%
%

%
%

\section{Supplementary: Saddle point method heuristics}
\label{sec:Heuristics}

A great reference for multivariate generating functions is Pemantle and Wilson \cite{PemantleWilson}.
Another fine reference for univariate generating functions is Flajolet and Sedgewick \cite{FlajoletSedgewick}.
We include the current analysis for the convenience of the reader.

\subsection{A  preliminary formal calculation}
Imagine that we have, due to Cauchy's integral formula, an equation of the form
$$
A_{\boldsymbol{\lambda}}\, =\, a_{\boldsymbol{\lambda}} \int_{-K_1(\boldsymbol{\lambda})}^{K_1(\boldsymbol{\lambda})}
\cdots \int_{-K_n(\boldsymbol{\lambda})}^{K_n(\boldsymbol{\lambda})} 
e^{-i \boldsymbol{q}_{\boldsymbol{\lambda}} \cdot \boldsymbol{x}} 
G_{\boldsymbol{\lambda}}(\boldsymbol{x})\, dx_1\, \cdots dx_n\, ,
$$
where $\boldsymbol{q}_{\boldsymbol{\lambda}}$ is a vector in $\R^n$ such that when $\|\boldsymbol{\lambda}\|$ is large.
then also $\boldsymbol{q}_{\boldsymbol{\lambda}}$ is large.

We assume that the generating function
satisfies that for each $\boldsymbol{x}$ we have
$$
G_{\boldsymbol{\lambda}}(\boldsymbol{x})
- \frac{1}{1-i \boldsymbol{q}_{\boldsymbol{\lambda}}\cdot \boldsymbol{x}+\frac{1}{2}\, \|\boldsymbol{x}\|^2} \to 0\, ,
$$
as $\boldsymbol{\lambda} \to \infty$
(where the ``convergence to infinity'' may need to be further delineated, such as we usually at least demand all components to ``go to 
infinity'' separately).
We also assume that $K_1(\boldsymbol{\lambda}),\dots,K_n(\boldsymbol{\lambda}) \to \infty$ as $\boldsymbol{\lambda} \to \infty$.

Note that the pointwise convergence that we have posited here is not sufficient for convergence of the integrals.
But in this subsection, we are proceeding formally (non-rigorously).
Then we expect that in some sense the integral is asymptotic to the formal integral (not obviously convergent):
$$
\int_{-K_1(\boldsymbol{\lambda})}^{K_1(\boldsymbol{\lambda})}
\cdots \int_{-K_n(\boldsymbol{\lambda})}^{K_n(\boldsymbol{\lambda})} 
e^{-i  \boldsymbol{q}_{\boldsymbol{\lambda}} \cdot \boldsymbol{x}} 
G_{\boldsymbol{\lambda}}(\boldsymbol{x})\, dx_1\, \cdots dx_n
\stackrel{?}{\sim}\, \text{\LARGE `` }
\int_{-\infty}^{\infty}
\cdots \int_{-\infty}^{\infty} 
\frac{ e^{-i  \boldsymbol{q}_{\boldsymbol{\lambda}} \cdot \boldsymbol{x}} } 
{1-i \boldsymbol{q}_{\boldsymbol{\lambda}}\cdot \boldsymbol{x}+\frac{1}{2}\, \|\boldsymbol{x}\|^2}\, 
dx_1\, \cdots dx_n\, . \text{\LARGE ''}
$$
A way to make sense of the latter is to use the ``Laplace transform trick,'' 
$$
\forall y \in \R\, ,\ \forall x \in (0,\infty)\, ,\ \text{ we have }\
\frac{1}{x+iy}\, =\, \int_0^{\infty} e^{-t\cdot (x+iy)}\, dt\, .
$$
Then if we interchange the order of integration (by a formal, unjustified step) we arrive at the desired integral
$$
\int_{-\infty}^{\infty} e^{-t} \left(\int_{-\infty}^{\infty}
\cdots \int_{-\infty}^{\infty} 
e^{-i  \boldsymbol{q}_{\boldsymbol{\lambda}} \cdot \boldsymbol{x}}
e^{-i t  \boldsymbol{q}_{\boldsymbol{\lambda}} \cdot \boldsymbol{x}} 
e^{-t\|\boldsymbol{x}\|^2/2}\, 
dx_1\, \cdots dx_n\right)\, dt\, .
$$
Having made many formal, non-rigorous steps, we could 
then complete the analysis by doing the Gaussian integrals to obtain
\begin{equation*}
\int_0^{\infty} e^{-t} \exp\left(-\frac{(1-t)^2}{2t}\, 
\| \boldsymbol{q}_{\boldsymbol{\lambda}}\|^2\right)\, \left(\frac{2\pi}{t}\right)^{n/2}\, dt\, 
\sim\, \frac{e^{-1} (2\pi)^{(n+1)/2}}{\| \boldsymbol{q}_{\boldsymbol{\lambda}}\|}\, ,\ 
\text{ as $\|\boldsymbol{\lambda}\|\to\infty$,}
\end{equation*}
where the final asymptotic equivalence follows by the usual Hayman admissibility analysis.
So, from all this we would expect 
\begin{equation}
\label{eq:oscSPM}
A_{\boldsymbol{\lambda}}\, \sim\, e^{-1} (2\pi)^{(n+1)/2}\, \cdot 
\frac{a_{\boldsymbol{\lambda}}}{\| \boldsymbol{q}_{\boldsymbol{\lambda}}\|}\, ,\ 
\text{ as $\|\boldsymbol{\lambda}\|\to\infty$.}
\end{equation}
\subsubsection{Example: binomial coefficient}
As an example, consider the binomial coefficient
$$
A_{(k,\ell)}\, =\, \binom{k+\ell}{k,\ell}\, ,
$$
which has the generating function 
\begin{equation*}
	M_2^{(1)}(x,y)\, =\, \sum_{k,\ell=0}^{\infty} \binom{k+\ell}{k,\ell} x^k y^{\ell}\, =\, \frac{1}{1-x-y}\, , 
\end{equation*}
using the nomenclature from (\ref{eq:M1formula}). By Cauchy's integral formula
\begin{equation*}
\begin{split}
\binom{k+\ell}{k,\ell}\, &=\, \oint_{\mathcal{C}(0;r)} \left(\oint_{\mathcal{C}(0;s)} \frac{1}{1-x-y} \cdot \frac{1}{x^k} \cdot \frac{1}{y^{\ell}}\, 
\frac{dy}{2\pi i y}\right)\, \frac{dx}{2\pi i x} \Bigg|_{\substack{r=k/(k+\ell+1)\\ s=\ell/(k+\ell+1)}}\\
&\qquad =\, \frac{(k+\ell+1)^{k+\ell+1}}{k^k\, \ell^{\ell}} \cdot \frac{1}{4\pi^2} 
\int_{-\pi}^{\pi} \int_{-\pi}^{\pi} 
\frac{e^{-i k \theta - i \ell \phi}}
{1-k(e^{i\theta}-1)-\ell(e^{i\phi}-1)}\, d\theta\, d\phi\, .
\end{split}
\end{equation*}
Now let us define $\boldsymbol{x} = (x_1,x_2) \in \R^2$ by the formula
$$
	(x_1,x_2)\, =\, (\sqrt{k}\, \theta\, ,\ \sqrt{\ell}\, \phi)\, ,\qquad 
(\theta,\phi)\, =\, (k^{-1/2}x_1,\ell^{-1/2}x_2)\, ,\qquad 
d\theta\, d\phi\, =\, \frac{1}{\sqrt{k\ell}}\, dx_1\, dx_2\, .
$$
Then we obtain
\begin{equation*}
\begin{split}
\binom{k+\ell}{k,\ell}\, &=\,  
\frac{(k+\ell+1)^{k+\ell+1}}{4\pi^2 k^{k+(1/2)}\, \ell^{\ell+(1/2)}} 
\int_{-\pi \sqrt{k}}^{\pi \sqrt{k}} \int_{-\pi \sqrt{\ell}}^{\pi \sqrt{\ell}} 
\frac{e^{-i \sqrt{k}\, x_1 - i \sqrt{\ell}\, x_2}}
{1-k(e^{ix_1/\sqrt{k}}-1)-\ell(e^{ix_2/\sqrt{\ell}}-1)}\, dx_1\, dx_2\, .
\end{split}
\end{equation*}
So in this case $\boldsymbol{\lambda} = (k,\ell)$ and $\boldsymbol{q}_{\boldsymbol{\lambda}}=(k^{1/2},\ell^{1/2})$.
We have $K_1(\boldsymbol{\lambda})=\pi\, \sqrt{k}$ and $K_2(\boldsymbol{\lambda})=\pi\, \sqrt{\ell}$.
We have
$$
a_{\boldsymbol{\lambda}}\, =\, \frac{(k+\ell+1)^{k+\ell+1}}{4\pi^2 k^{k+(1/2)}\, \ell^{\ell+(1/2)}}\,
=\, \left(1+\frac{1}{k+\ell}\right)^{k+\ell+1}\, \frac{(k+\ell)^{k+\ell+1}}{4\pi^2 k^{k+(1/2)}\, \ell^{\ell+(1/2)}}\,
\sim\, \frac{e}{4\pi^2}\, \cdot \frac{(k+\ell)^{k+\ell+1}}{k^{k+(1/2)}\, \ell^{\ell+(1/2)}}\, .
$$
Finally, by Taylor expansion, we have
\begin{equation*}
\begin{split}
G_{\boldsymbol{\lambda}}(\boldsymbol{x})\, 
&=\, \frac{1}{1-k(e^{ix_1/\sqrt{k}}-1)-\ell(e^{ix_2/\sqrt{\ell}}-1)}\\
&=\, \frac{1}{1-i\sqrt{k}\, x_1 - \frac{1}{2}\, x_1^2 -i \sqrt{\ell}\, x_2 - \frac{1}{2}\, x_2^2 +O(|x_1|^3/\sqrt{k})+O(|x_2|^3/\sqrt{\ell})}\, .
\end{split}
\end{equation*}
So it is true that pointwise $G_{\boldsymbol{\lambda}}(\boldsymbol{x})$ converges to 
$1/(1-i \boldsymbol{q}_{\boldsymbol{\lambda}}\cdot \boldsymbol{x}+\frac{1}{2}\, \|\boldsymbol{x}\|^2)$.
Note that $\|\boldsymbol{q}_{\boldsymbol{\lambda}}\|=\sqrt{k+\ell}$.
Thus, we are led to the purely formal guess that
$$
	A_{(k,\ell)}\, \sim\, e^{-1}\, (2\pi)^{3/2}\, \cdot \frac{e}{4\pi^2}\, \frac{(k+\ell)^{k+\ell+1}}{k^{k+(1/2)}\, \ell^{\ell+(1/2)}}\,
\cdot  \frac{1}{(k+\ell)^{1/2}}\, .
$$
This does give the correct answer of $(2\pi)^{-1/2} (k+\ell)^{k+\ell+(1/2)}/(k^{k+(1/2)}\ell^{\ell+(1/2)})$,
which is obtainable more easily using Stirling's formula.
The reason to consider this example, is that it is useful to keep it in mind, while taking on more complicated problems.

\section{Supplementary: Saddle point analysis  simplest assumptions}
\label{sec:Saddle}

The following assumptions are modeled on the multinomial coefficients.
A particularly simplifying assumption is that we have 
$$
G_{\boldsymbol{\lambda}}(\boldsymbol{x})\, =\, \frac{1}{g_{\boldsymbol{\lambda}}(\boldsymbol{x})}\, ,
$$
for a function $g_{\boldsymbol{\lambda}}(\boldsymbol{x})$ which is well-approximated pointwise by 
$1-i \boldsymbol{q}_{\boldsymbol{\lambda}}\cdot \boldsymbol{x}+\frac{1}{2}\, \|\boldsymbol{x}\|^2$,
when $\boldsymbol{\lambda} \to \infty$.
\begin{itemize}
\item[$\boldsymbol{\mathcal{A}1}$.]
We assume $G_{\boldsymbol{\lambda}}(\boldsymbol{x})\, =\, 1/g_{\boldsymbol{\lambda}}(\boldsymbol{x})$
where in turn
$g_{\boldsymbol{\lambda}}(\boldsymbol{x}) = 1 + \sum_{r=1}^{n} \gamma_r(x_r;\boldsymbol{\lambda})$
for functions $\gamma_1,\dots,\gamma_n$ satisfying
$$
\forall x_r \in [-K_1(\boldsymbol{\lambda}),K_1(\boldsymbol{\lambda})]\, ,\ \text{ we have }\
\operatorname{Re}[\gamma_r(x_r;\boldsymbol{\lambda})]\, \geq\, 0\, .
$$
Later we will actually assume a stronger condition, which we consider to be coercivity.
\end{itemize}
In particular this implies
$$
\inf\left(\left\{g_{\boldsymbol{\lambda}}(\boldsymbol{x})\, :\, 
\boldsymbol{x} \in [-K_1(\boldsymbol{\lambda}),K_1(\boldsymbol{\lambda})] \times \cdots \times [-K_n(\boldsymbol{\lambda}),K_n(\boldsymbol{\lambda})]\right\}\right)\, >\, 0\, .
$$
More precisely,  the infimum is actually $1$.
Then we may write
$$
A_{\boldsymbol{\lambda}}\, =\, a_{\boldsymbol{\lambda}} \int_{-K_1(\boldsymbol{\lambda})}^{K_1(\boldsymbol{\lambda})}
\cdots \int_{-K_n(\boldsymbol{\lambda})}^{K_n(\boldsymbol{\lambda})} 
e^{-i \boldsymbol{q}_{\boldsymbol{\lambda}} \cdot \boldsymbol{x}} 
\left(\int_0^{\infty} e^{-t g_{\boldsymbol{\lambda}}(\boldsymbol{x})}\, dt\right)\, dx_1\, \cdots dx_n\, ,
$$
using the Laplace transform trick. Moreover, by our assumption, we have that 
the function
$e^{-t g_{\boldsymbol{\lambda}}(\boldsymbol{x})}$
is absolutely summable on 
$$
(\boldsymbol{x},t) \in [-K_1(\boldsymbol{\lambda}),K_1(\boldsymbol{\lambda})] \times \cdots \times [-K_n(\boldsymbol{\lambda}),K_n(\boldsymbol{\lambda})] \times [0,\infty).
$$
Hence, we may interchange the integrals. 
So we have
$$
A_{\boldsymbol{\lambda}}\, =\, a_{\boldsymbol{\lambda}} 
\int_0^{\infty} \left(
\int_{-K_1(\boldsymbol{\lambda})}^{K_1(\boldsymbol{\lambda})}
\cdots \int_{-K_n(\boldsymbol{\lambda})}^{K_n(\boldsymbol{\lambda})} 
e^{-i \boldsymbol{q}_{\boldsymbol{\lambda}} \cdot \boldsymbol{x}} 
e^{-t g_{\boldsymbol{\lambda}}(\boldsymbol{x})}\,  dx_1\, \cdots dx_n\right)\, dt\, .
$$
Let write
the components of $\boldsymbol{q}_{\boldsymbol{\lambda}}$ in the following notation
\begin{equation}
\boldsymbol{q}_{\boldsymbol{\lambda}}\, =\, \big(q_1(\boldsymbol{\lambda}),\dots,q_n(\boldsymbol{\lambda})\big)\, .
\end{equation}
Now, rewriting $g_{\boldsymbol{\lambda}}(\boldsymbol{x}) = 1 + \sum_{r=1}^{n} \gamma_r(x_r;\boldsymbol{\lambda})$,
we have
$$
A_{\boldsymbol{\lambda}}\, =\, a_{\boldsymbol{\lambda}} 
\int_0^{\infty} e^{-t} \prod_{r=1}^{n}\left(
\int_{-K_r(\boldsymbol{\lambda})}^{K_r(\boldsymbol{\lambda})}
e^{-i q_r(\boldsymbol{\lambda}) x_r} 
e^{-t \gamma_r(x_r;\boldsymbol{\lambda})}\,  dx_r\right)\, dt\, ,
$$
where the $q_r(\boldsymbol{\lambda})$'s are the components of $\boldsymbol{q}_{\boldsymbol{\lambda}}$.
Let us denote
$$
\mathcal{I}_r(t;\boldsymbol{\lambda})\, =\, \int_{-K_r(\boldsymbol{\lambda})}^{K_r(\boldsymbol{\lambda})}
e^{-i q_r(\boldsymbol{\lambda}) x_r} 
e^{-t \gamma_r(x_r;\boldsymbol{\lambda})}\,  dx_r\, .
$$
So we have
$$
A_{\boldsymbol{\lambda}}\, =\, a_{\boldsymbol{\lambda}} 
\int_0^{\infty} e^{-t} \prod_{r=1}^{n} \mathcal{I}_r(t;\boldsymbol{\lambda})\, dt\, .
$$
\subsection{Digression to consider product nature of errors}
Recall that we are trying to prove that $\int_0^{\infty} e^{-t} \prod_{r=1}^{n} \mathcal{I}_r(t;\boldsymbol{\lambda})\, dt$
is asymptotic to $e^{-1} (2\pi)^{(n+1)/2}/\|\boldsymbol{q}_{\boldsymbol{\lambda}}\|$.
But our basic method will be to show that for a good family of $t$'s we have $\mathcal{I}_r(t;\boldsymbol{\lambda})$
is asymptotic to $\exp\left(-(1-t)^2q_r^2(\boldsymbol{\lambda})/(2t)\right)\left(2\pi/t\right)^{1/2}$, for each $r$.
But there will be errors. We may summarize this schematically by saying
$$
\mathcal{I}_r(t;\boldsymbol{\lambda})\, =\, e^{-(1-t)^2q_r^2(\boldsymbol{\lambda})/(2t)}\, \sqrt{\frac{2\pi}{t}}
+ \operatorname{Err}_r(t;\boldsymbol{\lambda})\, .
$$
The nature of the integral is such that as long as the errors are all uniformly bounded then we will be able to conclude
\begin{equation}
\begin{split}
\left|\int_0^{\infty} e^{-t} \left( \left(\prod_{r=1}^{n} \mathcal{I}_r(t;\boldsymbol{\lambda})\right)
-  \exp\left(-\frac{(1-t)^2}{2t}\, 
\| \boldsymbol{q}_{\boldsymbol{\lambda}}\|^2\right)\, \left(\frac{2\pi}{t}\right)^{n/2}\right)\, dt\right|\\
&\hspace{-10cm}
=\, \sum_{\substack{X \subset [n]\\ X\neq \emptyset, X \neq [n]}}
\int_0^{\infty} e^{-t}
\exp\left(-\frac{(1-t)^2}{2t}\, 
\sum_{r \in X} \boldsymbol{q}_{r}^2(\boldsymbol{\lambda})\right)\, \left(\frac{2\pi}{t}\right)^{|X|/2}\, 
\left(\prod_{s \in [n] \setminus X} \operatorname{Err}_s(t;\boldsymbol{\lambda})\right)\, dt\\
&\hspace{-8cm} + \int_0^{\infty} e^{-t} \left(\prod_{r=1}^n \operatorname{Err}_r(t;\boldsymbol{\lambda})\right)\, dt
\end{split}
\end{equation}
If we assume that for each $r \in [n]$ we have $\operatorname{Err}_r(t;\boldsymbol{\lambda}) \to 0$ pointwise for all $t \in [0,\infty)$
as $\boldsymbol{\lambda} \to \infty$, then we may do the integral to rewrite this as 
\begin{equation}
\label{eq:prodToReference}
\begin{split}
\left| \left(\int_0^{\infty} e^{-t} \left(\prod_{r=1}^{n} \mathcal{I}_r(t;\boldsymbol{\lambda})\right)
dt\right) - \frac{e^{-1} (2\pi)^{(n+1)/2}}{\|\boldsymbol{q}_{\boldsymbol{\lambda}}\|} \right|\,
&=\, \sum_{\substack{X \subset [n]\\ X\neq \emptyset, X \neq [n]}} o\left(\frac{1}
{\sqrt{\sum_{r \in X} \boldsymbol{q}_{r}^2(\boldsymbol{\lambda})}}\right)
\\
&\qquad + \int_0^{\infty} e^{-t} \left(\prod_{r=1}^n \operatorname{Err}_r(t;\boldsymbol{\lambda})\right)\, dt
\end{split}
\end{equation}
So that shows that essentially all we need to check going forward is that 
$$
\|\boldsymbol{q}_{\boldsymbol{\lambda}}\|\, \cdot \int_0^{\infty} e^{-t} \left|\prod_{r=1}^n \operatorname{Err}_r(t;\boldsymbol{\lambda})\right|\, dt\, 
\to\, 0\, ,\ \text{ as $\boldsymbol{\lambda} \to \infty$.}
$$
But it also motivates the next assumption.
\subsection{Resumption of statements of assumptions: small $t$ cutoff}

We are going to make the following uniformity assumption on the components of $\boldsymbol{q}_{\boldsymbol{\lambda}}$.
\begin{itemize}
\item[$\boldsymbol{\mathcal{A}2}$.] We assume that there exists a number $\eta \in (0,1)$ such that:
we have
$$
\eta\, <\, \frac{q_r(\boldsymbol{\lambda})}{q_s(\boldsymbol{\lambda})}\, <\, \frac{1}{\eta}\, ,
$$
for all $r,s \in \{1,\dots,n\}$, for all
$\boldsymbol{\lambda}$'s in the range we consider. 
\end{itemize}
This way we can think of the components of $\boldsymbol{q}_{\boldsymbol{\lambda}}$ as each being as large as the norm
$\|\boldsymbol{q}_{\boldsymbol{\lambda}}\|$ modulo a constant factor.
In particular then we can rewrite equation (\ref{eq:prodToReference}) as 
\begin{equation*}
\begin{split}
\left| \left(\int_0^{\infty} e^{-t} \left(\prod_{r=1}^{n} \mathcal{I}_r(t;\boldsymbol{\lambda})\right)
dt\right) - \frac{e^{-1} (2\pi)^{(n+1)/2}}{\|\boldsymbol{q}_{\boldsymbol{\lambda}}\|} \right|\,
&=\, o\left(\frac{1}{\|\boldsymbol{q}_{\boldsymbol{\lambda}}\|}\right)
+ \int_0^{\infty} e^{-t} \left(\prod_{r=1}^n \operatorname{Err}_r(t;\boldsymbol{\lambda})\right)\, dt
\end{split}
\end{equation*}

Now we are goint to partition $t$ into a small $t$ regime, a large $t$ regime and an intermediate regime.
We make this assumption.
\begin{itemize}
\item[$\boldsymbol{\mathcal{A}3}$.] For each $p \in (0,1)$, we denote a small $t$ cutoff as
$$
\tau_p(\boldsymbol{\lambda})\, =\, \|\boldsymbol{q}_{\boldsymbol{\lambda}}\|^{-p}\, .
$$
We assume that we have 
$$
\|\boldsymbol{q}_{\boldsymbol{\lambda}}\|\int_0^{\tau_p(\boldsymbol{\lambda})} e^{-t}\, \prod_{r=1}^{n} |\mathcal{I}_r(t;\boldsymbol{\lambda})|\, dt\, \to\, 0\ \text{ as $\boldsymbol{\lambda}\to\infty$,}
$$
for each $p \in (0,1)$, separately.
\end{itemize}
What will occur in examples is that we will have an ad-hoc argument -- based on analyticity
of $G_{\boldsymbol{\lambda}}$ and the Euler-Maclaurin summation formula -- to show that 
for each $m>1$ we have
$|\mathcal{I}_r(t;\boldsymbol{\lambda})|$ is bounded by  a constant times $(t+|q_r(\boldsymbol{\lambda})|^{-1})^m$.
Using $\tau_p(\boldsymbol{\lambda})\, =\, \|\boldsymbol{q}_{\boldsymbol{\lambda}}\|^{-p}$,
if we take $m$ high enough then we obtain the desired bound.

Let us delay further description until the relevant section. Since there is rapid oscillation of the factor
$e^{-i q_r(\boldsymbol{\lambda}) x_r}$ in the integral
$\mathcal{I}_r(t;\boldsymbol{\lambda})\, dt$, which is $\int_{-K_r(\boldsymbol{\lambda})}^{K_r(\boldsymbol{\lambda})}
e^{-i q_r(\boldsymbol{\lambda}) x_r} 
e^{-t \gamma_r(x_r;\boldsymbol{\lambda})}\,  dx_r$,
it might not be surprising that when $t$ is small the total integral is small.
For example, if $t=0$, then the integral of just the factor $e^{-i q_r(\boldsymbol{\lambda}) x_r}$
is $0$ because it will be the case that $[-K_r(\boldsymbol{\lambda}),K_r(\boldsymbol{\lambda})]$
comprises exactly an integer number of full periods.
(Another intuitive hint is the fact that the main term of the integral is $e^{-(1-t)^2q_r^2(\boldsymbol{\lambda})/(2t)}\, \sqrt{\frac{2\pi}{t}}$
which converges to $0$ as $t \to 0^+$ faster than any power.)

\subsection{Large $|x_r|$ regime ($t$ not too small)}

In examples, the $\gamma_r$ functions will turn out to be coercive. 
\begin{itemize}
\item[$\boldsymbol{\mathcal{A}4}$.] There is a constant $\sigma^2\in [1,\infty)$ such that 
$$
\operatorname{Re}[\gamma_r(x_r;\boldsymbol{\lambda})]\, \geq\, \frac{x_r^2}{2\sigma^2}\, .
$$
\end{itemize}
Recall that we are trying to establish that the integral
$$
\mathcal{I}_r(t;\boldsymbol{\lambda})\,  =
\int_{-K_r(\boldsymbol{\lambda})}^{K_r(\boldsymbol{\lambda})}
e^{-i q_r(\boldsymbol{\lambda}) x_r} 
e^{-t \gamma_r(x_r;\boldsymbol{\lambda})}\,  dx_r\, 
$$
is well approximated by $e^{-(1-t)^2 q_r^2(\boldsymbol{\lambda})/(2t)}\, \sqrt{2\pi/t}$.
By the coercivity, for each $\rho \in (0,1)$, we will be able to find a $M_{\rho}(\boldsymbol{\lambda})$ which is  polylogarithmic in $\|\boldsymbol{q}_{\boldsymbol{\lambda}}\|$
such that 
$$
\|\boldsymbol{q}_{\boldsymbol{\lambda}}\|^{1+\rho}\, (\sqrt{t}\, M_{\rho}(\boldsymbol{\lambda}))^n
\prod_{r=1}^{n}\left(
\int_{[-K_r(\boldsymbol{\lambda}),K_r(\boldsymbol{\lambda})]\setminus[-M_{\rho}(\boldsymbol{\lambda})/\sqrt{t},M_{\rho}(\boldsymbol{\lambda})/\sqrt{t}]}
|e^{-t \gamma_r(x_r;\boldsymbol{\lambda})}|\,  dx_r\right)\, \to\, 0\, , \text{ as $\boldsymbol{\lambda}\to\infty$.}
$$
The reason for the prefactor $(\sqrt{t}\, M_{\rho}(\boldsymbol{\lambda}))^n$ is that a standard bound gives 
$$
\int_a^{\infty} e^{-x^2/2}\, dx\, \leq\, \frac{e^{-a^2/2}}{a}\, ,\ \text{ for $a>0$.}
$$
But then by scaling
$$
\int_a^{\infty} e^{-tx^2/2}\, dx\, \leq\, \frac{e^{-ta^2/2}}{at}\, ,\ \text{ for $a>0$.}
$$
But we are choosing $a=M_{\rho}(\boldsymbol{\lambda})/\sqrt{t}$ in order to have $e^{-nta^2/(2\sigma^2)} \leq \|\boldsymbol{q}_{\boldsymbol{\lambda}}\|^{-(1+\rho)}$ (in which case we may take 
$M_{\rho}(\boldsymbol{\lambda})=\sqrt{2(1+\rho)\sigma^2 \ln(\|\boldsymbol{q}_{\boldsymbol{\lambda}}\|)/n}$).
Then the prefactor $(\sqrt{t}\, M_{\rho}(\boldsymbol{\lambda}))^n$ is due to the $n$th power of the denominator
$at=\sqrt{t}\, M_{\rho}(\boldsymbol{\lambda})$.

However, the number $M_{\rho}(\boldsymbol{\lambda})$ is large, not small. Also, $t\geq \tau_p(\boldsymbol{\lambda})$ in the intermediate
zone. And this is $\|\boldsymbol{q}_{\boldsymbol{\lambda}}\|^{-p}$ for a $p$ as small as desired.
Therefore, choosing $\rho>pn/2$ (assuming $p$ is quite close to $0$), we still have that for $t$ not too small
$$
\|\boldsymbol{q}_{\boldsymbol{\lambda}}\|\, 
\prod_{r=1}^{n}\left(
\int_{[-K_r(\boldsymbol{\lambda}),K_r(\boldsymbol{\lambda})]\setminus[-M_{\rho}(\boldsymbol{\lambda})/\sqrt{t},M_{\rho}(\boldsymbol{\lambda})/\sqrt{t}]}
|e^{-t \gamma_r(x_r;\boldsymbol{\lambda})}|\,  dx_r\right)\, \to\, 0\, , \text{ as $\boldsymbol{\lambda}\to\infty$.}
$$
\subsection{Intermediate regime for $x_r$ (and $t$ not too small)}
Now we turn attention to $t$ in the  regime $[\tau_p(\boldsymbol{\lambda}),\infty)$
and also $x_r$ in the intermediate $x$-regime $[-M_{\rho}(\boldsymbol{\lambda})/\sqrt{t},M_{\rho}(\boldsymbol{\lambda})/\sqrt{t}]$.
Then, we assume that the following holds.
\begin{itemize}
\item[$\boldsymbol{\mathcal{A}5}$.] We assume there is a constant $c$ such that
$$
\gamma_r(x_r;\boldsymbol{\lambda}) + i q_r(\boldsymbol{\lambda}) x_r - \frac{1}{2}\, x_r^2\, 
=\, \frac{c x_r^3}{q_r(\boldsymbol{\lambda})}
+ O\left(\frac{M_{\rho}^4(\boldsymbol{\lambda})}{t^{2}q_r^2(\boldsymbol{\lambda})}\right)\, ,
$$
for all $t \in [\tau_p(\boldsymbol{\lambda}),T_{\boldsymbol{\lambda}}]$ and 
$x \in [-M_{\rho}(\boldsymbol{\lambda})/\sqrt{t},M_{\rho}(\boldsymbol{\lambda})/\sqrt{t}]$.
\end{itemize}
This assumption will follow from Taylor expansion in our examples.
By the statements above, $t$ is no smaller than an arbitrarily small power of $1/\|\boldsymbol{q}_{\boldsymbol{\lambda}}\|$ 
and $M_{\rho}(\boldsymbol{\lambda})$  is polylogarithmic in $\|\boldsymbol{q}_{\boldsymbol{\lambda}}\|$.
So, the error which is the final term on the right-hand-side of the formula above should be small. 
Let us denote the error as 
$$
\mathcal{E}_r(x_r;\boldsymbol{\lambda})\, =\, 
\gamma_r(x_r;\boldsymbol{\lambda}) + i q_r(\boldsymbol{\lambda}) x_r - \frac{1}{2}\, x_r^2 - \frac{c x_r^3}{q_r(\boldsymbol{\lambda})}\, .
$$
So Assumption 6 insures that
$$
\mathcal{E}_r(x_r;\boldsymbol{\lambda})\, =\, 
O\left(\frac{M_{\rho}^4(\boldsymbol{\lambda})}{t^{2}q_r^2(\boldsymbol{\lambda})}\right)\, .
$$
But now we can rewrite the integrand in $\mathcal{I}_r(t;\boldsymbol{\lambda})$ as 
$$
e^{-i q_r(\boldsymbol{\lambda}) x_r} 
e^{-t \gamma_r(x_r;\boldsymbol{\lambda})}\, =\, 
e^{-i (1-t)q_r(\boldsymbol{\lambda}) x_r
-\frac{1}{2}\, tx_r^2} e^{-tcx_r^3q_r^{-1}(\boldsymbol{\lambda})-t\mathcal{E}_r(x_r;\boldsymbol{\lambda})}\, .
$$
And actually both extra terms are small in the intermediate $t$ and $x_r$ regime: $tcx_r^3q_r^{-1}(\boldsymbol{\lambda})$
and $t\mathcal{E}_r(x_r;\boldsymbol{\lambda})$.
Note that in the intermediate regime $x_r = O(M_{\rho}(\boldsymbol{\lambda})/\sqrt{t})$ so that
$tcx_r^3q_r^{-1}(\boldsymbol{\lambda})=O\left(\frac{M_{\rho}^3(\boldsymbol{\lambda})}
{q_r(\boldsymbol{\lambda})\, \sqrt{t}}\right)$.
So we may Taylor expand the final exponential as 
$$
e^{-tcx_r^3q_r^{-1}(\boldsymbol{\lambda})-t\mathcal{E}_r(x_r;\boldsymbol{\lambda})}\,
=\, 1 -\frac{ctx_r^3}{q_r(\boldsymbol{\lambda})} + O\left(\frac{M_{\rho}^4(\boldsymbol{\lambda})}{tq_r^2(\boldsymbol{\lambda})}\right)\, .
$$
The term containing the constant $c$ just barely fails to be small enough to neglect in the special case that $n=1$
(which for the case of multinomial coefficients is the case of a sequence of $1$'s).
But by keeping it and doing the exact Gaussian integral we can see, {\em a posteriori}, that it does not contribute a significant source of error.
(In the case $n\geq 2$ then it would be squared or multiplied to a higher power, in which case, it could be neglected even before taking the Gaussian integral.)

Using coercivity outside the intermediate $x_r$-regime, as in the last subsection, we will then obtain
\begin{equation}
\label{eq:IncludesGaussian}
\int_{-M_{\rho}(\boldsymbol{\lambda})/\sqrt{t}}^{M_{\rho}(\boldsymbol{\lambda})/\sqrt{t}}
e^{-i q_r(\boldsymbol{\lambda}) x_r} 
e^{-t \gamma_r(x_r;\boldsymbol{\lambda})}\,  dx_r\, 
=\, \int_{-\infty}^{\infty} e^{-i (1-t)q_r(\boldsymbol{\lambda}) x_r-\frac{1}{2}\, tx_r^2}\left(1 -\frac{ctx_r^3}{q_r(\boldsymbol{\lambda})}\right)\, dx_r + \widetilde{\operatorname{Err}}_r(t;\boldsymbol{\lambda})\, ,
\end{equation}
where
$$
|\widetilde{\operatorname{Err}}_r(t;\boldsymbol{\lambda})|\,
\leq\, \left(1+\frac{|c|M_{\rho}^3(\boldsymbol{\lambda})}{t^{3/2}q_r(\boldsymbol{\lambda})}\right) 
\widehat{\operatorname{Err}}_r(t;\boldsymbol{\lambda})
+O\left(\frac{M_{\rho}^3(\boldsymbol{\lambda})}{t^{3/2} q_r^2(\boldsymbol{\lambda})}\right)
$$
where the nonnegative term $\widehat{\operatorname{Err}}_r(t;\boldsymbol{\lambda})$ is an error as in the last subsection, satisfying
$$
\|\boldsymbol{q}_{\boldsymbol{\lambda}}\|^{1+\rho}\, (\sqrt{t}\, M_{\rho}(\boldsymbol{\lambda}))^n
\prod_{r=1}^{n}\left(\widehat{\operatorname{Err}}_r(t;\boldsymbol{\lambda})\right)\, \to\, 0\, , \text{ as $\boldsymbol{\lambda}\to\infty$.}
$$
This is good enough to ignore all of $\widetilde{\operatorname{Err}}_r(t;\boldsymbol{\lambda})$.
So, having dispensed with the errors other than the exact Gaussian integral, we return to the integral on the right-hand-side of 
(\ref{eq:IncludesGaussian}):
\begin{equation*}
\begin{split}
\int_{-\infty}^{\infty} e^{-i (1-t)q_r(\boldsymbol{\lambda}) x_r-\frac{1}{2}\, tx_r^2}\left(1 -\frac{ctx_r^3}{q_r(\boldsymbol{\lambda})}\right)\, dx_r\\ 
&\hspace{-4cm}=\, 
e^{-(1-t)^2q_r^2(\boldsymbol{\lambda})/(2t)}\, \sqrt{\frac{2\pi}{t}}\left(1 + \frac{ct}{q_r(\boldsymbol{\lambda})} \cdot \frac{c_1 (1-t)q_r(\boldsymbol{\lambda})}{\sqrt{t}}
\left(\frac{c_2 (1-t)^2 q_r^2(\boldsymbol{\lambda})+c_3}{t}\right)\right)\, .
\end{split}
\end{equation*}
But in the usual Hayman analysis, for any powers $m_1,\dots,m_n \in \{0,1,\dots\}$, we have
$$
\int_0^{\infty} e^{-t} \exp\left(-\frac{(1-t)^2}{2t}\, 
\| \boldsymbol{q}_{\boldsymbol{\lambda}}\|^2\right)\, \left(\frac{2\pi}{t}\right)^{n/2}\, 
\prod_{r=1}^{n} \left(\frac{(1-t)q_r(\boldsymbol{\lambda})}{\sqrt{t}}\right)^{m_r} dt\, 
= O\left(\frac{1}{\| \boldsymbol{q}_{\boldsymbol{\lambda}}\|}\right)\, ,
$$
in part using Assumption ${\mathcal{A}2}$.
Therefore, because of the extra factor of $ct/q_r(\boldsymbol{\lambda})$,
we may ignore the error term arising from the cubic power $x_r^3$ in the Gaussian integral.

This takes care of all the technical details, modulo checking the technical assumptions.

\section{Supplementary: Multinomial coefficients toy problem}
\label{sec:Multinomial}

In order to illustrate the saddle-point method in a multivariate setting,  let us consider 
$$
A_{\boldsymbol{\lambda}}\, =\, \binom{\lambda_1+\dots+\lambda_n}{\lambda_1,\dots,\lambda_n}\, ,
$$
for $\boldsymbol{\lambda} \in \{0,1,\dots\}^n$.
In order to satisfy $\mathcal{A}2$, let us assume that there are numbers $\Lambda_1,\dots,\Lambda_n \in (0,1]$
such that the mode in which $\boldsymbol{\lambda} \to \infty$ is such that
$$
	\lambda_1,\dots,\lambda_n \to \infty\, ,\ \text{ such that }\ \frac{\lambda_r}{\lambda_1+\dots+\lambda_r} \to \Lambda_r\, ,\
\text{ for $r \in [n]$.}
$$
The generating function is one of the most elementary generating functions 
\begin{equation}
	M_n^{(1)}(\boldsymbol{z})\, =\, \sum_{\boldsymbol{\lambda} \in \{0,1,\dots\}^n}
\binom{\lambda_1+\dots+\lambda_n}{\lambda_1,\dots,\lambda_n} z_1^{\lambda_1} \cdots z_n^{\lambda_n}\, ,
\end{equation}
using the nomenclature from (\ref{eq:M1formula}).
So
$$
	M_n^{(1)}(\boldsymbol{z})\, =\, \frac{1}{1-(z_1+\dots+z_n)}\, .
$$
So, using Cauchy's integral formula and choosing the radii of the circles such that the radius for $z_r$ is $\lambda_r/(1+\lambda_1+\dots+\lambda_n)$, we have
$$
\binom{\lambda_1+\dots+\lambda_n}{\lambda_1,\dots,\lambda_n}\,
=\, \alpha_{\boldsymbol{\lambda}}\, 
\int_{-\pi}^{\pi} \cdots \int_{-\pi}^{\pi} 
\frac{e^{-i \boldsymbol{\lambda} \cdot \boldsymbol{\theta}}}{1-\sum_{r=1}^{n} \lambda_r (e^{i\theta_r}-1)}\, 
d\theta_1 \cdots d\theta_n\, ,
$$
where
$$
\alpha_{\boldsymbol{\lambda}}\, =\, 
(2\pi)^{-n}\, \frac{(1+\lambda_1+\dots+\lambda_n)^{1+\lambda_1+\dots+\lambda_n}}{\lambda_1^{\lambda_1} \cdots \lambda_n^{\lambda_n}}\, .
$$
Then rescaling, to let $x_r = \theta_r\, \sqrt{\lambda_r}$ for each $r \in [n]$, then
$$
A_{\boldsymbol{\lambda}}\, =\, a_{\boldsymbol{\lambda}} \int_{-K_1(\boldsymbol{\lambda})}^{K_1(\boldsymbol{\lambda})}
\cdots \int_{-K_n(\boldsymbol{\lambda})}^{K_n(\boldsymbol{\lambda})} 
e^{-i \boldsymbol{q}_{\boldsymbol{\lambda}} \cdot \boldsymbol{x}} 
\left(\int_0^{\infty} e^{-t g_{\boldsymbol{\lambda}}(\boldsymbol{x})}\, dt\right)\, dx_1\, \cdots dx_n\, ,
$$
where
$$
a_{\boldsymbol{\lambda}}\, =\, \frac{\alpha_{\boldsymbol{\lambda}}}{\sqrt{\lambda_1\cdots\lambda_n}}\, =\,
(2\pi)^{-n}\, \frac{(1+\lambda_1+\dots+\lambda_n)^{1+\lambda_1+\dots+\lambda_n}}{\lambda_1^{\lambda_1+(1/2)} \cdots \lambda_n^{\lambda_n+(1/2)}}\, ,
$$
$q_r(\boldsymbol{\lambda}) = \sqrt{\lambda_r}$, $K_r(\boldsymbol{\lambda}) = \pi\, \sqrt{\lambda_r}$
and
$$
g_{\boldsymbol{\lambda}}(\boldsymbol{x})\, =\, 1 + \sum_{r=1}^{n} \gamma_r(x_r;\boldsymbol{\lambda})\, ,
$$
for
$$
\gamma_r(x_r;\boldsymbol{\lambda})\, =\, -\lambda_r \left(e^{i x_r/\sqrt{\lambda_r}}-1\right)\, .
$$
So Assumption $\mathcal{A}1$ is satisfied.
Note that elementary trigonometry leads to 
$$
\gamma_r(x_r;\boldsymbol{\lambda})\, =\, 2\lambda_r \sin^2\left(\frac{x_r}{2\sqrt{\lambda_r}}\right)
-i \lambda_r\, \sin\left(\frac{x_r}{\sqrt{\lambda_r}}\right)\, .
$$
So we have the coercivity condition
\begin{equation}
\label{eq:MultiCoerc}
\operatorname{Re}[\gamma_r(x_r;\boldsymbol{\lambda})]\, =\, 2\lambda_r \sin^2\left(\frac{x_r}{2\sqrt{\lambda_r}}\right)\,
\geq\, \frac{2x_r^2}{\pi^2}\, ,
\end{equation}
on the interval $x_r \in [-\pi\, \sqrt{\lambda_r}\, ,\ \pi\, \sqrt{\lambda_r}]$ by convexity/concavity of $\sin(\theta)$ for $\theta \in [0,\pi/2]$.
So condition $\mathcal{A}4$ is satisfied with $\sigma^2 = \pi^2/4$.
Hence, condition $\mathcal{A}1$ is also satisfied (because it is actually subsumed by condition $\mathcal{A}4$).
We may also note that condition $\mathcal{A}5$ is also satisfied because
$$
\gamma_r(x_r;\boldsymbol{\lambda})\, =\, -\lambda_r \left(e^{i x_r/\sqrt{\lambda_r}}-1\right)\, =\,
-i x_r\, \sqrt{\lambda_r} + \frac{1}{2}\, x_r^2 + \frac{i x_r^3}{6\, \sqrt{\lambda_r}} + O\left(\frac{x_r^4}{\lambda_r}\right)\, .
$$
So the constant in $\mathcal{A}5$ is $c=i/6$.

\subsection{Checking the small $t$ condition: analyticity and Euler-Maclaurin}
Now let us consider the hardest condition to check, assumption $\mathcal{A}3$.
Here we explain the ad-hoc argument using analyticity and the Euler-Maclaurin summation formula.
Firstly, note that the range of integration $x_r \in [-\pi\, \sqrt{\lambda_r}\, ,\ \pi\, \sqrt{\lambda_r}]$ is an integer multiple of the period of 
$e^{-i\sqrt{\lambda_r}\, x_r}$. Of course this was no coincidence since we just rescaled an integral originating in Cauchy's integral formula.
Now let us consider the integral over one period. We are motivated by the idea that for small $t$, the function 
$e^{-t \gamma_r(x_r;\boldsymbol{\lambda})}$ may not vary much. (And a constant time $e^{-i\sqrt{\lambda_r}\, x_r}$ integrates
to $0$ over exactly 1 period.) So we have
$$
\int_{2\pi n/\sqrt{\lambda_r}}^{2\pi (n+1)/\sqrt{\lambda_r}} e^{-i x_r\, \sqrt{\lambda_r}} 
e^{-t\gamma_r(x_r;\boldsymbol{\lambda})}\, dx_r\,
=\, \int_{2\pi n/\sqrt{\lambda_r}}^{2\pi (n+1)/\sqrt{\lambda_r}} e^{-i x_r\, \sqrt{\lambda_r}} 
\exp\left(-t\lambda_r\, (e^{ix_r/\sqrt{\lambda_r}}-1)\right)\, dx_r\, .
$$
Now we note that the general $m$th order derivative formula of the second factor 
(for $m \in \{1,2,\dots\}$) may be written
$$
\frac{d^m}{dx_r^m} \exp\left(-t\sqrt{\lambda_r}\, (e^{ix_r/\sqrt{\lambda_r}}-1)\right)\,
=\, \sum_{\mu=1}^{m} c_{m,\mu}\,  t^{\mu} \lambda_r^{\mu-(m/2)}\,  \exp\left(-t\sqrt{\lambda_r}\, (e^{ix_r/\sqrt{\lambda_r}}-1)\right) e^{i\mu x_r/\sqrt{\lambda_r}}\, ,
$$
for some constants $c_{m,\mu}$.
So, denoting by $X^{(r)}_n = 2 \pi n/ \sqrt{\lambda_r}$, we may use Taylor's theorem to $M$th order to write
\begin{equation*}
\begin{split}
e^{-t\gamma_r(x_r;\boldsymbol{\lambda})} -  e^{-t\gamma_r(X_n^{(r)};\boldsymbol{\lambda})}\, 
&=\,
e^{-t\gamma_r(X_n^{(r)};\boldsymbol{\lambda})}
\sum_{m=1}^{M} \sum_{\mu=1}^{m} \frac{c_{m,\mu}}{m!}\, \cdot  t^{\mu} \lambda_r^{\mu-(m/2)}\, 
(x_r-X^{(r)}_n)^m  \,  
e^{i\mu X^{(r)}_n/\sqrt{\lambda_r}}\\
&\hspace{-1cm} + O\left(\max_{X^{(r)}_n\leq x_r \leq X^{(r)}_{n+1}}|e^{-t\gamma_r(x_r;\boldsymbol{\lambda})}|
t(t+\lambda_r^{-1})^{M}\lambda_r^{(M+1)/2}(X^{(r)}_{n+1}-X^{(r)}_n)^{M+1}\right)\, ,
\end{split}
\end{equation*}
for $X^{(r)}_n\leq x_r \leq X^{(r)}_{n+1}$.
Since $X^{(r)}_{n+1}-X^{(r)}_{n}=2\pi/\sqrt{\lambda_r}$, we may rewrite this as 
\begin{equation*}
\begin{split}
e^{-t\gamma_r(x_r;\boldsymbol{\lambda})}-e^{-t\gamma_r(X_n^{(r)};\boldsymbol{\lambda})}\, 
&=\, e^{-t\gamma_r(X_n^{(r)};\boldsymbol{\lambda})}
\sum_{m=1}^{M} \sum_{\mu=1}^{m} \frac{c_{m,\mu}}{m!}\, \cdot  t^{\mu} \lambda_r^{\mu-(m/2)}\, 
(x_r-X^{(r)}_n)^m \,  
e^{i\mu X^{(r)}_n/\sqrt{\lambda_r}}\\
&\qquad\qquad + O\left(\max_{X^{(r)}_n\leq x_r \leq X^{(r)}_{n+1}}|e^{-t\gamma_r(x_r;\boldsymbol{\lambda})}|
t(t+\lambda_r^{-1})^{M}\right)\, ,
\end{split}
\end{equation*}
Then using coercivity (\ref{eq:MultiCoerc}) to resum the final error term as a Riemann sum, we may write
\begin{equation*}
\begin{split}
\mathcal{I}_r(t;\boldsymbol{\lambda})\,  &=\, O\left((t+\lambda_r^{-1})^{M}\, \sqrt{t}\right)\\
&\hspace{-2cm}
+\sum_{\mu=1}^{M} \sum_{n=-N_{r,1}}^{N_{r,2}}
 e^{-t\gamma_r(X_n^{(r)};\boldsymbol{\lambda})}
e^{i\mu X^{(r)}_n/\sqrt{\lambda_r}}
\sum_{m=\mu}^{M} \frac{c_{m,\mu}}{m!}\, \cdot  t^{\mu} \lambda_r^{\mu-(m/2)}\, 
\int_{X_n^{(r)}}^{X_{n+1}^{(r)}} (x_r-X^{(r)}_n)^m  e^{-ix_r\, \sqrt{\lambda_r}}\, dx_r \, ,
\end{split}
\end{equation*}
where
$N_{r,1} = \lfloor\lambda_r/2\rfloor$ and 
$N_{r,2}=\lceil \lambda_r/2 \rceil - 1$.
We left out the constant times $e^{ix_r\, \sqrt{\lambda_r}}$ integral because it evaluates to $0$.
Performing the integrals, we may write
\begin{equation*}
\begin{split}
\mathcal{I}_r(t;\boldsymbol{\lambda})\,  &=\, \frac{2\pi}{\sqrt{\lambda_r}}\, 
\sum_{\mu=1}^{M} \sum_{n=-N_{r,1}}^{N_{r,2}}
 e^{-t\gamma_r(X_n^{(r)};\boldsymbol{\lambda})}
e^{i\mu X^{(r)}_n/\sqrt{\lambda_r}}
\sum_{m=\mu}^{M}  \frac{c_{m,\mu}}{m!}\, \cdot  t^{\mu} \lambda_r^{-(m-\mu)}J_m
+O\left((t+\lambda_r^{-1})^{M}\, \sqrt{t}\right)\, , 
\end{split}
\end{equation*}
where $J_m = \int_0^{2\pi} \theta^m e^{-i\theta}\, d\theta$.
We view the remaining summation $\frac{2\pi}{\sqrt{\lambda_r}} \sum_{n=-N_{r,1}}^{N_{r,2}}$ as a Riemann sum approximation to a Riemann integral.
It remains to show that this integral is no larger than a high power of $t$.
But that follows by analyticity and the Euler-Maclaurin summation formula.

More specifically, by analyticity, we note that
$$
\int_{-\pi\sqrt{\lambda_r}}^{\pi \sqrt{\lambda_r}} e^{-t\gamma_r(x_r;\boldsymbol{\lambda})}
e^{i\mu x_r/\sqrt{\lambda_r}}\, dx\, =\, 0\, ,
$$
since the source of a pole, the term $e^{-ix_r\, \sqrt{\lambda_r}}$, is now gone (integrated out).
So the question is how well does the Riemann sum approximate the Riemann integral.
The typical approach to this question is to use the Euler-Maclaurin summation formula.
The endpoint values of $e^{-t\gamma(x_r;\boldsymbol{\lambda})}$, evaluated at the upper limit of summation
minus the lower limit of summation, will always cancel due to periodicity.
So then, the Euler-Maclaurin summation formula may be written, using 
$$
f_{r,\mu}(x_r)\, =\, e^{-t\gamma_r(x_r;\boldsymbol{\lambda})}
e^{i\mu x_r/\sqrt{\lambda_r}}\, ,
$$
we have (due to the fact that the comparison Riemann integral equals 0) just the remainder term:
$$
\frac{2\pi}{\sqrt{\lambda_r}}\, 
\sum_{n=-N_{r,1}}^{N_{r,2}}f_{r,\mu}(X^{(r)}_n)\,
=\, (-1)^{p+1} \int_{-\pi\, \sqrt{\lambda}}^{\pi\, \sqrt{\lambda}} f^{(p)}_{r,\mu}(x_r)\, 
\frac{P_p(x_r\, \sqrt{\lambda_r}/(2\pi))}{p!}\, dx_r\, ,
$$
where $P_p$ is the periodized version of the $p$th Bernoulli polynomial. In particular, it is bounded
by a constant just depending on $p$.
But taking the $p$th derivative of $f_{r,\mu}$ is similar to taking the $p$th derivative of $f_{r,0}$
which we did at the outset of this subsection. We find:
$$
f_{r,\mu}(x_r)^{(p)}\, =\, O(t(t+\lambda_r^{-1})^{p-1} e^{-t\gamma_r(x_r;\boldsymbol{\lambda}_r)})\, .
$$
So, again using coercivity (\ref{eq:MultiCoerc}) we have
$$
\frac{2\pi}{\sqrt{\lambda_r}}\, 
\sum_{n=-N_{r,1}}^{N_{r,2}}f_{r,\mu}(X^{(r)}_n)\,
=\, O\left((t+\lambda_r^{-1})^{p-1}\, \sqrt{t}\right)\, .
$$
Putting this together with the previous bound, we get
$$
\mathcal{I}_r(t;\boldsymbol{\lambda})\,  =\, 
O\left((t+\lambda_r^{-1})^{M}\, \sqrt{t}\right) 
+ O\left((t+\lambda_r^{-1})^{p-1}\, \sqrt{t}\right)\, .
$$
This is what we desired. We may take $t \leq \tau_p(\boldsymbol{\lambda})=\|\boldsymbol{\lambda}\|^{-p}$ and sufficiently high powers of $M$ and $p$, to obtain what we want.  

\subsection{Conclusion: verification of the well known formula}

So we have by equation (\ref{eq:oscSPM}) and the formula for $a_{\boldsymbol{\lambda}}$:
$$
A_{\boldsymbol{\lambda}}\, \sim\, e^{-1} (2\pi)^{(n+1)/2}\, \cdot 
\frac{a_{\boldsymbol{\lambda}}}{\| \boldsymbol{q}_{\boldsymbol{\lambda}}\|}\, \sim\,
e^{-1} (2\pi)^{(n+1)/2}\, \cdot
\frac{1}{\sqrt{\lambda_1+\dots+\lambda_n}}\, \cdot
(2\pi)^{-n}\, \frac{(1+\lambda_1+\dots+\lambda_n)^{1+\lambda_1+\dots+\lambda_n}}{\lambda_1^{\lambda_1+\frac{1}{2}} \cdots \lambda_n^{\lambda_n+\frac{1}{2}}}\, .
$$
But as before, we may write $(1+\lambda_1+\dots+\lambda_n)^{1+\lambda_1+\dots+\lambda_n}$
as $(\lambda_1+\dots+\lambda_n)^{1+\lambda_1+\dots+\lambda_n}$ times 
$$
\left(1+\frac{1}{\lambda_1+\dots+\lambda_n}\right)^{1+\lambda_1+\dots+\lambda_n}\, \sim\, e\, ,\ \text{ as $\boldsymbol{\lambda} \to \infty$.}
$$
So we have
$$
A_{\boldsymbol{\lambda}}\,
=\, \binom{\lambda_1+\dots+\lambda_n}{\lambda_1,\dots,\lambda_n}\,
\sim\, (2\pi)^{-(n-1)/2}\, \cdot 
\frac{(\lambda_1+\dots+\lambda_n)^{\lambda_1+\dots+\lambda_n+\frac{1}{2}}}{\lambda_1^{\lambda_1+\frac{1}{2}} \cdots \lambda_n^{\lambda_n+\frac{1}{2}}}\,  ,\ \text{ as $\boldsymbol{\lambda} \to \infty$,}
$$
where recall that our mode of convergence is such that $\lambda_1,\dots,\lambda_n \to \infty$ but 
$\lambda_r/(\lambda_1+\dots+\lambda_n) \to \Lambda_r$ for numbers $\Lambda_1,\dots,\Lambda_n \in (0,1]$.

\section{Auxiliary: Proof of LD rate function for $\mathcal{A}(\kappa N,\kappa N,\gamma N)$}
Here we prove equation (\ref{eq:SymmetricAAsymptotics}). Let us replace $K$ by $k_N$, and let us replace $j_K$ by $j_N$.
Then with that notation, we are trying to prove the following version of the previous formula: assuming $K_N \sim \kappa N$
and $j_N \sim \gamma N$ as $N \to \infty$, for $\kappa,\gamma \in (0,\infty)$, we have
\begin{equation}
\label{eq:SymmetricA2}
\lim_{N \to \infty} \frac{1}{2K_N}\, \ln\left(\mathcal{A}(K_N,K_N,j_N)\right)\, 
=\, 2 \ln(2) + (1+\rho) \ln(1+\rho) - \rho \ln(\rho)\Bigg|_{\rho=\gamma/(4\kappa)}\, .
\end{equation}
\subsection{Setup using the Cauchy integral formula }
By Corollary \ref{cor:GF}, we know that the generating function for $\mathcal{A}(k,\ell,j)$
is $((1-2(x+y)+(x-y)^2)^{1/2} - z)$.
By the Cauchy integral formula, we know that 
$$
\mathcal{A}(k,k,j)\, =\, \int_{-\pi}^{\pi} s^{-j} e^{-ij\phi}
\left(\int_{-\pi}^{\pi} \int_{-\pi}^{\pi} \frac{t^{-2k} e^{-ik(\theta_1+\theta_2)}}
{\mathbbm{D}(te^{i\theta_1},t e^{i\theta_2},se^{i\phi})}\, \frac{d\theta_1}{2\pi}\, \frac{d\theta_2}{2\pi}\right)\,
\frac{d\phi}{2\pi}\, ,
$$
for 
$$
{\mathbbm{D}(te^{i\theta_1},t e^{i\theta_2},se^{i\phi})}\, =\, 
\left(1-2t(e^{i\theta_1}+e^{i\theta_2})
+t^2(e^{i\theta_1}-e^{i\theta_2})^2\right)^{1/2} - s e^{i\phi}\, .
$$
This is true
for any $t>0$ and $s>0$ as long as: $t<1/4$ and $s<(1-4t)^{1/2}$.

Recall that we are taking $k=k_N \sim \kappa N$ and $j=j_N \sim \gamma N$ as $N \to \infty$.
Then our choice is 
$$
t\, =\, t_N\, =\, \frac{\kappa}{4\kappa+\gamma} - \frac{\kappa}{(4\kappa+\gamma)^2}\, N^{-1} + o(N^{-1})\,
$$
and
$$
s\, =\, s_N\, =\, \frac{\sqrt{\gamma}}{\sqrt{4\kappa+\gamma}} 
- \frac{2\kappa+\gamma}{\sqrt{\gamma}\, (4\kappa+\gamma)^{3/2}}\, N^{-1} + o(N^{-1})\, .
$$
In a later section we will explain how to arrive at these formulas. The saddle point method leads to an optimization
problem for $t$ and $s$, as is usual.
But the solution of that equation is more algebraic than analytical.
So we delay its explanation.

Let us define
$$
t_*\, =\, \frac{\kappa}{4\kappa+\gamma}\ \text{ and }\ 
s_*\, =\, \sqrt{\frac{\gamma}{4\kappa+\gamma}}\, .
$$
Then $t_N<t_*$, which implies $(1-4t_N)^{1/2}>(1-4t_*)^{1/2}$.
And $s_N<s_*$. Hence, $t=t_N$ and $s=s_N$ are valid choices.
Of course, 
$$
t_*\, =\, \lim_{N \to \infty} t_N\ \text{ and }\ s_*\, =\, \lim_{N \to \infty} s_N\, .
$$
Then we note that, the leading-order part of the Cauchy integral formula is
\begin{equation}
\label{eq:LDR1}
\lim_{N \to \infty} \frac{1}{N}\, \ln\left(\frac{t_N^{-2k_N} s_N^{-j_N}}{\mathbbm{D}(t_N,t_N,s_N)}\right)\,
=\, -2\kappa \ln(t_*) - \gamma \ln(s_*)\, ,
\end{equation}
because $\mathbb{D}(t_N,t_N,s_N)$ only behaves algebraically in $N$, not exponentially.
Let us explain that by first rationalizing the denominator (which we need for later):
$$
\frac{1}{\mathbb{D}(z_1,z_2,w)}\, =\, \frac{\mathcal{N}(z_1,z_2,w)}{\mathcal{D}(z_1,z_2,w)}\, ,
$$
for
$\mathcal{N}(z_1,z_2,w) = (1-2(z_1+z_2)+(z_1-z_2)^2)^{1/2} + w$
and
$\mathcal{D}(z_1,z_2,w) = 1-2(z_1+z_2)+(z_1-z_2)^2-w^2$.
Then
$$
\mathcal{N}(t_*,t_*,s_*)\, =\, (1-4t_*)^{1/2} + s_*\, =\, 2\, \sqrt{\frac{\gamma}{4\kappa+\gamma}}\, ,
$$
while
$$
\mathcal{D}(t_N,t_N,s_N)\, =\, 1-4t_N - s_N^2\, \sim\, \frac{2}{4\kappa+\gamma}\, \cdot \frac{1}{N}\, .
$$
So
$$
\frac{1}{\mathbb{D}(t_N,t_N,s_N)}\, \sim\, \sqrt{\gamma(4\kappa+\gamma)}\, N\, ,
$$
showing the validity of (\ref{eq:LDR1}).
But then that means
$$
\lim_{N \to \infty} \frac{1}{N}\, \ln\left(\frac{t_N^{-2k_N} s_N^{-j_N}}{\mathbbm{D}(t_N,t_N,s_N)}\right)\,
=\, -2\kappa \Big(\ln(\kappa)-\ln(4\kappa+\gamma)\Big) - \gamma \Big(\frac{1}{2}\, \ln(\gamma) - \frac{1}{2}\, \ln(4\kappa+\gamma)\Big)\, .
$$
If we use $\gamma = 4 \rho \kappa$ this gives
$$
\lim_{N \to \infty} \frac{1}{N}\, \ln\left(\frac{t_N^{-2k_N} s_N^{-j_N}}{\mathbbm{D}(t_N,t_N,s_N)}\right)\,
=\, -2\kappa \Big(\ln(\kappa)-\ln(4\kappa+4\kappa\rho)\Big) 
- 4\kappa \rho \Big(\frac{1}{2}\, \ln(4\kappa\rho) - \frac{1}{2}\, \ln(4\kappa+4\kappa\rho)\Big)\, .
$$
Simplifying this formula gives the right-hand-side of (\ref{eq:SymmetricA2}).
But note that we are not done.
We are  just beginning, because we must now prove that the remaining integral,
$$
\int_{-\pi}^{\pi} 
\int_{-\pi}^{\pi} \int_{-\pi}^{\pi} \frac{e^{-ij\phi-ik(\theta_1+\theta_2)}}
{\mathbbm{D}(t_Ne^{i\theta_1},t_N e^{i\theta_2},s_Ne^{i\phi})/\mathbbm{D}(t_N,t_N,s_N)}\, 
\frac{d\theta_1}{2\pi}\, \frac{d\theta_2}{2\pi}\,
\frac{d\phi}{2\pi}\, ,
$$
is neither exponentially large nor exponentially small in $N$.
We do this following the general methodology of  Sections \ref{sec:Heuristics}, \ref{sec:Saddle}, \ref{sec:Multinomial}.

\subsection{Checking coercivity after combining roots on the unit circle}

There are two things to check, coercivity and the small Laplace parameter behavior which we previously called small $t$ bounds
but which could also be called ``infrared bounds.'' The small Laplace parameter behavior is via an ad hoc argument related to analyticity
and the Euler-Maclaurin summation formula. But technically we should establish coercivity first, in order to make sure that we can
use the Laplace transform trick.

Coercivity is usually easy to check by trigonometry methods. This time, because we rationalized the denominator,
we introduced another root for $\mathcal{D}(t_*,t_*,s_*e^{i\phi})$ at $\phi=\pi$.
So we want to combine those two roots; thereby, we are reducing the range of $\phi$ from $[-\pi,\pi)$ to $[-\pi/2,\pi/2)$.
Let us define
$$
\mathcal{N}_0(z_1,z_2,w)\, =\, \big(1-2(z_1+z_2)+(z_1-z_2)\big)^{1/2}\qquad \text{ and }\qquad
\mathcal{N}_1(z_1,z_2,w)\, =\, w\, ,
$$
so that $\mathcal{N}(z_1,z_2,w) = \mathcal{N}_0(z_1,z_2,w) + \mathcal{N}_1(z_1,z_2,w)$. Then we see that
$$
\mathcal{D}(z_1,z_2,-w)\, =\, \mathcal{D}(z_1,z_2,w)\qquad \text{ and }\qquad
\mathcal{N}(z_1,z_2,-w)\, =\, \mathcal{N}_0(z_1,z_2,w) - \mathcal{N}_1(z_1,z_2,w)\, .
$$
So, defining $d(j) = 0$ if $j$ is even and $d(j)=1$ if $j$ is odd (the remainder after division by 2) we have
\begin{equation}
\label{eq:Restart1}
\begin{split}
\int_{-\pi}^{\pi} 
\int_{-\pi}^{\pi} \int_{-\pi}^{\pi} \frac{e^{-ij\phi-ik(\theta_1+\theta_2)}}
{\mathbbm{D}(t_Ne^{i\theta_1},t_N e^{i\theta_2},s_Ne^{i\phi})}\, 
\frac{d\theta_1}{2\pi}\, \frac{d\theta_2}{2\pi}\,
\frac{d\phi}{2\pi}\\
&\hspace{-5cm}
=\,
\int_{-\pi/2}^{\pi/2} 
\left(\int_{-\pi}^{\pi} \int_{-\pi}^{\pi} \frac{e^{-ij\phi-ik(\theta_1+\theta_2)} \mathcal{N}_{d(j)}(t_Ne^{i\theta_1},t_Ne^{i\theta_2},s_Ne^{i\phi})}
{\mathcal{D}(t_Ne^{i\theta_1},t_N e^{i\theta_2},s_Ne^{i\phi})}\, 
\frac{d\theta_1}{2\pi}\, \frac{d\theta_2}{2\pi}\right)\,
\frac{d\phi}{\pi}\, .
\end{split}
\end{equation}
The two parts of the numerator $\mathcal{N}_0(z_1,z_2,w)$ and $\mathcal{N}_1(z_1,z_2,w)$
are equal at the saddle point $(t_N,t_N,s_N)$. That is the nature of rationalizing-the-denominator.
Their leading order behavior is $\sqrt{1-4t_*}=s_*=\sqrt{\gamma/(4\kappa+\gamma)}$, which is order-1.
It is neither vanishing nor diverging.
Therefore, the alternation between $\mathcal{N}_0(z_1,z_2,w)$ and $\mathcal{N}_1(z_1,z_2,w)$
does not much affect the asymptotics of $\mathcal{A}(k,k,j)$. But at the next order, it would. (We could start to see ``number theoretic''
effects in the asymptotics at higher order than the leading large-deviation-rate-function behavior.)

Note that on both sides of equation (\ref{eq:Restart1}) we have eliminated a factor of $\mathbbm{D}(t_N,t_N,s_N)$ compared to what
we had been considering prior.
That is because such a factor is just algebraic in $N$.
We are just trying to check that the integral is neither exponentially large nor exponentially small in $N$.
So such factors do not concern us at this stage of the argument.
Therefore, we put-in another factor to consider the integral
$$
\int_{-\pi/2}^{\pi/2} 
\left(\int_{-\pi}^{\pi} \int_{-\pi}^{\pi} \frac{e^{-ij\phi-ik(\theta_1+\theta_2)} \mathcal{N}_{d(j)}(t_Ne^{i\theta_1},t_Ne^{i\theta_2},s_Ne^{i\phi})}
{\mathcal{D}(t_Ne^{i\theta_1},t_N e^{i\theta_2},s_Ne^{i\phi})/\mathcal{D}(t_N,t_N,s_N)}\, 
\frac{d\theta_1}{2\pi}\, \frac{d\theta_2}{2\pi}\right)\,
\frac{d\phi}{\pi}\, .
$$
So, at this stage, we are just going to check coercivity of 
$$
\frac{\mathcal{D}(t_Ne^{i\theta_1},t_N e^{i\theta_2},s_Ne^{i\phi})}{\mathcal{D}(t_N,t_N,s_N)}
$$
for $(\theta_1,\theta_2,\phi) \in [-\pi,\pi]\times [-\pi,\pi] \times [-\pi/2,\pi/2]$.
We note that the fraction can be written as 
$$
\frac{\mathcal{D}(t_Ne^{i\theta_1},t_N e^{i\theta_2},s_Ne^{i\phi})}{\mathcal{D}(t_N,t_N,s_N)}\,
=\, 1-a_N(e^{i\theta_1}+e^{i\theta_2}-2)
+b_N(e^{i\theta_1}-e^{i\theta_2})^2 - c_N (e^{2i\phi}-1)\, ,
$$
for $a_N = 2t_N/\mathcal{D}(t_N,t_N,s_N)$, $b_N = t_N^2/\mathcal{D}(t_N,t_N,s_N)$ and $c_N=s_N/\mathcal{D}(t_N,t_N,s_N)$.
These satisfy
$$
a_N\, \sim\, \kappa N\, ,\qquad 
b_N\, \sim\,  \frac{\kappa^2}{2(4\kappa+\gamma)}\, N\, ,\qquad
c_N\, \sim\, \frac{\gamma}{2}\, N\, ,\qquad \text{ as $N \to \infty$.}
$$
Let us denote
$$
\mathfrak{d}(e^{i\theta_1},e^{i\theta_2},e^{i\phi})\,
=\, a_N\, \left(e^{i\theta_1}+e^{i\theta_2}-2\right) -b_N\, \left(e^{i\theta_1}-e^{i\theta_2}\right)^2
+ c_N\, \left(e^{2i\phi}-1\right)\, .
$$
By elementary trigonometry, this can be written
$$
\mathfrak{d}(e^{i\theta_1},e^{i\theta_2},e^{i\phi})\,
=\, \mathfrak{d}_{\Re}(\theta_1,\theta_2,\phi) + i \mathfrak{d}_{\Im}(\theta_1,\theta_2,\phi)\, ,
$$
for
\begin{equation}
\begin{split}
\mathfrak{d}_{\Re}(\theta_1,\theta_2,\phi)\, 
&=\, 
2 a_N\, \left(\sin^2\left(\frac{\theta_1}{2}\right) +\sin^2\left(\frac{\theta_2}{2}\right) \right)\\
&\qquad  - 4 b_N\, \cos\left(\frac{\theta_1+\theta_2}{2}\right) \sin^2\left(\frac{\theta_1-\theta_2}{2}\right)\\
&\qquad + 2c_N\, \sin^2\left(\phi\right)\, ,
\end{split}
\end{equation}
and
\begin{equation}
\begin{split}
\mathfrak{d}_{\Im}(\theta_1,\theta_2,\phi)\, 
&=\, 
a_N\, \left(\sin\left(\theta_1\right) +\sin\left(\theta_2\right) \right)\\
&\qquad  + 4 b_N\, \cos\left(\frac{\theta_1+\theta_2}{2}\right) \sin^2\left(\frac{\theta_1-\theta_2}{2}\right)\\
&\qquad + c_N\, \sin\left(2\phi\right)\, .
\end{split}
\end{equation}
Now to check coercivity, note that $\sin^2(\alpha-\beta)\leq 2\sin^2(\alpha)+2\sin^2(\beta)$. Hence,
$$
\mathfrak{d}_{\Re}(\theta_1,\theta_2,\phi)\, \geq\, (2a_N-8b_N) \left(\sin^2\left(\frac{\theta_1}{2}\right)+\sin^2\left(\frac{\theta_2}{2}\right)\right) + 2 c_N\, \sin^2(\phi)\, ,
$$
which is positive because $2a_N-8b_N \sim \frac{4\kappa(2\kappa+\gamma)}{2(4\kappa+\gamma)}\, N$.
So, in turn
\begin{equation}
\label{eq:coerc1}
\mathfrak{d}_{\Re}(\theta_1,\theta_2,\phi)\, \geq\, (2a_N-8b_N) \left(\frac{\theta_1^2}{\pi^2}+\frac{\theta_2^2}{\pi^2}\right)
+ 8 c_N\, \frac{\phi^2}{\pi^2}\, ,
\end{equation}
for $(\theta_1,\theta_2,\phi) \in [-\pi,\pi)\times [-\pi,\pi) \times \left[-\frac{\pi}{2}\, ,\ \frac{\pi}{2}\right)$.

Now suppose that we take $(y_1,y_2,y_3) = (\theta_1\, \sqrt{N}\, ,\ \theta_2\, \sqrt{N}\, ,\ \theta_3\, \sqrt{N})$.
Then we have the integral
$$
\int_{0}^{\infty} e^{-\lambda} \mathcal{I}_N(\lambda)\, d\lambda\, ,
$$
where, in turn, $\mathcal{I}_N(\lambda)$ is the result of the iterated integral
$$
\int_{-L_N^{(1)}}^{L_N^{(1)}}
\int_{-L_N^{(2)}}^{L_N^{(2)}}
\int_{-L_N^{(3)}}^{L_N^{(3)}}
e^{-i\, \boldsymbol{p}_N \cdot \boldsymbol{y}}
\widetilde{\mathcal{N}}_{N,d(j)}(\boldsymbol{y})\, e^{-\lambda \widetilde{\mathcal{D}}_N(\boldsymbol{y})}\, dy_1\, dy_2\, dy_3
$$
where $L_N^{(1)} = L_N^{(2)} = \pi\sqrt{N}$ and $L_N^{(3)} = \frac{\pi}{2}\, \sqrt{N}$.
Here $\boldsymbol{p}_N = (k_N/\sqrt{N}\, ,\ k_N/\sqrt{N}\, ,\ j_N/\sqrt{N})$. And
$$
\widetilde{\mathcal{N}}_{N,d(j)}(\boldsymbol{y})\, =\, \mathcal{N}_{d(j)}\left(\frac{\boldsymbol{y}}{\sqrt{N}}\right)\ \text{ and }\
\widetilde{\mathcal{D}}_N(\boldsymbol{y})\, =\, \left(\mathcal{D}\left(\frac{\boldsymbol{y}}{\sqrt{N}}\right)/\mathcal{D}(t_N,t_N,s_N)\right)-1\, .
$$
In particular, from equation (\ref{eq:coerc1})
\begin{equation*}
\operatorname{Re}\left[\widetilde{\mathcal{D}}_N(\boldsymbol{y})\right]\, \geq\, \frac{2a_N-8b_N}{\pi^2 N}\,  \left(y_1^2+y_2^2\right)
+ \frac{8 c_N}{\pi^2 N}\, y_3^2\, 
\sim\, \frac{2\kappa (2\kappa+\gamma)}{(4\kappa+\gamma)\pi^2}\, (y_1^2+y_2^2) + 
\frac{4\gamma}{\pi^2}\, y_3^2\ \text{ as $N \to \infty$,}
\end{equation*}
for $(y_1,y_2,y_3) \in [-L_N^{(1)},L_N^{(1)}]\times [-L_N^{(2)},L_N^{(2)}] \times [-L_N^{(3)},L_N^{(3)}]$.

\subsection{Short digression: choice of coordinates and their transformation}
Note that the coercivity bound is not the actual quadratic behavior near the origin. 
The latter is what contributes to the leading order component of the final value of the integral.
The Taylor expansion out to quadratic order is as follows:
\begin{equation*}
\begin{split}
\widetilde{\mathcal{D}}_N(\boldsymbol{y})\, 
&=\, -i\left(\frac{a_N}{\sqrt{N}}\, \left(y_1+y_2\right)+\frac{2c_N}{\sqrt{N}}\, y_3\right) \\
&\qquad + \frac{a_N}{2N}\, \left(y_1^2+y_2^2\right) - \frac{b_N}{N}\, \left(y_1-y_2\right)^2 + \frac{c_N}{N}\, y_3^2 
+ O\left(\frac{1}{\sqrt{N}}\right)
\end{split}
\end{equation*}
which can also be written as 
\begin{equation*}
\widetilde{\mathcal{D}}_N(\boldsymbol{y})\, 
=\, -i\, \sqrt{N}\, (\kappa y_1 + \kappa y_2 + \gamma y_3)  + \frac{\kappa}{2}\, \left(y_1^2+y_2^2\right) - \frac{\kappa^2}{2(4\kappa+\gamma)}\, \left(y_1-y_2\right)^2 
+ \frac{\gamma}{2}\, y_3^2 
+ O\left(\frac{1}{\sqrt{N}}\right)\, ,
\end{equation*}
where we note that also
$\boldsymbol{p}_N \cdot \boldsymbol{y} \sim \sqrt{N}\, (\kappa y_1 + \kappa y_2 + \gamma y_3)$,
as $N \to \infty$.
At this point one might wish to rescale the coordinates $y_1,y_2,y_3$, which can be done.
But one might also wish to rotate 45-degrees in the $(y_1,y_2)$ plane in order to bring the quadratic component
into diagonal form.
That should not be done before dispensing with the bounds on $\mathcal{I}_N(\lambda)$ for small $\lambda$.

A rotation would complicate the form of the domain of integration $[-L_N^{(1)},L_N^{(1)}]\times [-L_N^{(2)},L_N^{(2)}] \times [-L_N^{(3)},L_N^{(3)}]$, which is  of product form. This is important because the limiting form of the integrand is not absolutely summable, although
of course the finite integrand on the domain of integration is absolutely summable. 
As opposed to the Hayman method for Gaussian type saddle points, we cannot merely cut-off the domain at a sphere
or ellipsoid. For the small Laplace transform parameter bounds that we need to establish, the structure of the domain is essential.
To carry out the analysis
we need the domain of integration to be an exact integer number of fundamental domains for the complex exponential.
This allows us to get cancellation.

\subsection{Bounds for small values of the Laplace parameter}

We believe that analyticity plus the Euler-Maclaurin summation formula gives the correct ad-hoc argument for the small $\lambda$
behavior of $\mathcal{I}_N(\lambda)$: namely, for the bounds. But now that we are in multi-dimensions, we need a multi-dimensional
Euler-Maclaurin summation formula. Undoubtedly this exists in the literature. In a future paper we will seek this out and use it.
But for now, let us note that an alternative perspective is to use the Fourier transform picture.
Then the Geometric sums identity (for a complex exponential sum) and elementary (Sobolev-type) inequalities
showing rapid decay of Fourier coefficients of an analytic function  can allow a bypass.
The underlying reasoning is the same: if we change to a sum of fundamental domains picture to integrate out the rapid oscillation
of $\exp(-i\boldsymbol{p}_N \cdot \boldsymbol{y})$, then after integrating out we get a sum that resembles a Riemann sum of 
an analytic function over a product of closed contours in its domain of analyticity.
So, comparing the sum to the integral, we reduce to a small error.

We remind ourselves that the reason we need small Laplace transform parameter bounds is that the coercivity which was established
in the last subsection is multiplied by $\lambda$ in the exponential. So at small values of $\lambda$ its strength becomes diluted.
However, the rapid oscillation of $\exp(-i\boldsymbol{p}_N \cdot \boldsymbol{y})$ essentially kicks-in everywhere except at $\lambda=1$.
At $\lambda=1$ the equal and opposite rapid oscillation of the leading part of $\exp(-\lambda \widetilde{\mathcal{D}}_N(\boldsymbol{y}))$,
which is $i\lambda \boldsymbol{p}_N \cdot \boldsymbol{y}$, cancels this.
In particular, we use that to control the small behavior for $0\leq \lambda\ll 1$, as $N \to \infty$.

\subsubsection{Sum of fundamental domains}

Recall that
we are trying to bound $\mathcal{I}_N(\lambda)$ for small $\lambda$, where
$$
\mathcal{I}_N(\lambda)\, =\, \int_{-L_N^{(1)}}^{L_N^{(1)}}
\int_{-L_N^{(2)}}^{L_N^{(2)}}
\int_{-L_N^{(3)}}^{L_N^{(3)}}
e^{-i\, \boldsymbol{p}_N \cdot \boldsymbol{y}}
\widetilde{\mathcal{N}}_{N,d(j)}(\boldsymbol{y})\, e^{-\lambda \widetilde{\mathcal{D}}_N(\boldsymbol{y})}\, dy_1\, dy_2\, dy_3\, ,
$$
where $L_N^{(1)} = L_N^{(2)} = \pi\sqrt{N}$ and $L_N^{(3)} = \frac{\pi}{2}\, \sqrt{N}$.
Let $\mathcal{V}_n =  \{-\lfloor n/2 \rfloor\, ,\ -\lfloor n/2 \rfloor+1\, ,\ \dots\, ,\ \lceil n/2\rceil - 1\}$.
Then, for $\boldsymbol{\nu} = (\nu_1,\nu_2,\nu_3) \in \mathcal{V}_k\times \mathcal{V}_k \times \mathcal{V}_j$, let us define
$$
\widehat{\mathcal{I}}_N(\lambda;\boldsymbol{\nu})\, =\, 
\iiint_{\mathcal{R}_{N,k,j}(\boldsymbol{\nu})}
e^{-i(ky_1+ky_2+jy_3)/\sqrt{N}}
\widetilde{\mathcal{N}}_{N,d(j)}(\boldsymbol{y}) e^{-\lambda \widetilde{\mathcal{D}}_N(\boldsymbol{y})}\, dy_1\, dy_2\, dy_3\, .
$$
where
$$
\mathcal{R}_{N,k,j}(\boldsymbol{\nu})\, =\, \left[\frac{2\pi\nu_1\sqrt{N}}{k},\frac{2\pi (\nu_1+1)\sqrt{N}}{k}\right]\times
\left[\frac{2\pi\nu_2\sqrt{N}}{k},\frac{2\pi (\nu_2+1)\sqrt{N}}{k}\right]\times
\left[\frac{2\pi\nu_3\sqrt{N}}{j},\frac{2\pi (\nu_3+1)\sqrt{N}}{j}\right]\, .
$$
\subsubsection{Taylor expansion}
Let us rescale this to be an order-1 domain integral. We do this so that the Taylor expansion is more transparent. 
Since $k\sim \kappa N$ and $j\sim \gamma N$, that means we could rescale by $\sqrt{N}$. So
$$
\widehat{\mathcal{I}}_N(\lambda;\boldsymbol{\nu})\, =\, 
N^{-3/2} \iiint_{\sqrt{N}\, \mathcal{R}_{N,k,j}(\boldsymbol{\nu})}
e^{-i(ky_1+ky_2+jy_3)/N}
\widetilde{\mathcal{N}}_{N,d(j)}\left(\frac{\boldsymbol{y}}{\sqrt{N}}\right) 
e^{-\lambda \widetilde{\mathcal{D}}_N(\boldsymbol{y}/\boldsymbol{\sqrt{N}})}\, dy_1\, dy_2\, dy_3\, ,
$$
where $\sqrt{N}\, \mathcal{R}_{N,k,j}(\boldsymbol{\nu}) = \{\sqrt{N}\, \boldsymbol{y}\, :\, \boldsymbol{y} \in \mathcal{R}_{N,k,j}(\boldsymbol{\nu})\}$, which is 
$$
\left[\frac{2\pi\nu_1N}{k},\frac{2\pi (\nu_1+1)N}{k}\right]\times
\left[\frac{2\pi\nu_2N}{k},\frac{2\pi (\nu_2+1)N}{k}\right]\times
\left[\frac{2\pi\nu_3N}{j},\frac{2\pi (\nu_3+1)N}{j}\right]\, .
$$
So this has order-1 3d volume.
But now, if we take a Taylor expansion about the vertex
$$
\boldsymbol{v}_{N,k,j}(\boldsymbol{\nu})\, =\, \left(\frac{2\pi\nu_1N}{k}\, ,\ \frac{2\pi\nu_3 N}{k}\, ,\ \frac{2\pi \nu_3 N}{j}\right)\, ,
$$ 
out to order-$n$, we get
\begin{equation*}
\begin{split}
\widetilde{\mathcal{N}}_{N,d(j)}\left(\frac{\boldsymbol{y}}{\sqrt{N}}\right) 
e^{-\lambda \widetilde{\mathcal{D}}_N(\boldsymbol{y}/\boldsymbol{\sqrt{N}})}\\
&\hspace{-4cm}
=\, \sum_{\substack{\boldsymbol{a} \in \{0,1,\dots\}^3\, :\\
|\boldsymbol{a}| \leq n}} \mathcal{M}_{N,k,j}^{({\boldsymbol{a}})}(\boldsymbol{\nu}) \left(y_1-\frac{2\pi\nu_1N}{k}\right)^{a_1} 
\left(y_2-\frac{2\pi\nu_2N}{k}\right)^{a_2} 
\left(y_3-\frac{2\pi\nu_3N}{j}\right)^{a_3}\\
&\hspace{-2cm}
+ \mathcal{E}_{N,k,j}^{(n)}(\boldsymbol{\nu},\boldsymbol{y})\, ,
\end{split}
\end{equation*}
where $|\boldsymbol{a}| = a_1+a_2+a_3$ and where the coefficients $\mathcal{M}_{N,k,j}^{({\boldsymbol{a}})}$ are obtained from the derivatives
$$
\mathcal{M}_{N,k,j}^{({\boldsymbol{a}})}\, =\, \frac{1}{(a_1!)(a_2!)(a_3!)}\, \cdot
\frac{\partial^{a_1+a_2+a_3}}{\partial y_1^{a_1} \, \partial y_2^{a_2} \, \partial  y_3^{a_3}} \left(\widetilde{\mathcal{N}}_{N,d(j)}\left(\frac{\boldsymbol{y}}{\sqrt{N}}\right) 
e^{-\lambda \widetilde{\mathcal{D}}_N(\boldsymbol{y}/\boldsymbol{\sqrt{N}})}\right)\Bigg|_{\boldsymbol{y} = \boldsymbol{v}_{N,k,j}(\boldsymbol{\nu})}\, .
$$
But now we note that we may rewrite this, using the original definitions of $\widetilde{\mathcal{N}}_N$ and $\widetilde{\mathcal{D}}_N$
$$
\mathcal{M}_{N,k,j}^{({\boldsymbol{a}})}\, =\, \frac{1}{(a_1!)(a_2!)(a_3!)}\, \cdot
\frac{\partial^{a_1+a_2+a_3}}{\partial y_1^{a_1} \, \partial y_2^{a_2} \, \partial  y_3^{a_3}} \left({\mathcal{N}}_{d(j)}\left(\frac{\boldsymbol{y}}{N}\right) 
e^{-\lambda \mathcal{D}(t_N,t_N,s_N)^{-1} {\mathcal{D}}(\boldsymbol{y}/\boldsymbol{N})}\right)\Bigg|_{\boldsymbol{y} = \boldsymbol{v}_{N,k,j}(\boldsymbol{\nu})}\, .
$$
So it is easy to see that $\mathcal{M}_{N,k,j}^{({\boldsymbol{a}})}$ is $O(\lambda^{|\boldsymbol{a}|}+N^{-|\boldsymbol{a}|})$.
In particular, the final error is of the size:
$$
\mathcal{E}_{N,k,j}^{(n)}(\boldsymbol{\nu},\boldsymbol{y})\, =\, O\left(\lambda^{n+1}+N^{-(n+1)}\right)\, .
$$
We are going to let the small $\lambda$ cut-off be of the form $\Lambda_{N}^{(p)}=N^{-p}$ for some small but fixed $p \in (0,1)$. 
So, since we will only apply the above argument
for $\lambda \in [0,\Lambda_N^{(p)}]$, the error is of the order $O(N^{-p(n+1)})$ which is small enough if we take $n$ to be large enough,
depending on the fixed $p$.
When we integrate the monomials against the complex exponentials, we get constants.
For the constant monomial $\boldsymbol{a}=(0,0,0)$ the integral is $0$.
\subsubsection{Consideration of the Fourier components and elementary bounds}
So we have 
$$
\widehat{\mathcal{I}}_N(\lambda;\boldsymbol{\nu})\, =\, O(N^{-\frac{3}{2} + (n+1)p}) + 
N^{-3/2} 
\sum_{\substack{\boldsymbol{a} \in \{0,1,\dots\}^3\, :\\
0<|\boldsymbol{a}| \leq n}} 
\Gamma_{N,k,j}^{(\boldsymbol{a})} \cdot 
\mathcal{M}_{N,k,j}^{({\boldsymbol{a}})}(\boldsymbol{\nu})\, ,
$$
where the $\Gamma_{N,k,j}^{(\boldsymbol{a})}$ constants are the result of integrating the monomial $y_1^{a_1} y_2^{a_2} y_3^{a_3}$
against $e^{-i(ky_1+ky_2+jy_3)/N}$ over the domain $[0,2\pi N/k] \times [0,2\pi N/k] \times [0,2\pi N/j]$.
But now we come to the main point of ``analyticity'' and the sum versus the integral.
The constants $\mathcal{M}_{N,k,j}^{({\boldsymbol{a}})}(\boldsymbol{\nu})$ are actually functions of $\boldsymbol{\nu}$.
They are defined by taking certain partial derivatives of the product of analytic functions 
$$
 \left({\mathcal{N}}_{d(j)}\left(\frac{\boldsymbol{y}}{N}\right) 
e^{-\lambda \mathcal{D}(t_N,t_N,s_N)^{-1} {\mathcal{D}}(\boldsymbol{y}/\boldsymbol{N})}\right)\, .
$$
As such, they are expressed in terms of Fourier modes.
But the only Fourier modes which are present in $\mathcal{M}_{N,k,j}^{({\boldsymbol{a}})}(\boldsymbol{\nu})$
have all nonnegative components, and must have at least 1 positive component.
There are no negative Fourier modes,
because the functions $\mathcal{N}$ and $\mathcal{D}$ are analytic.
Moreover, the zero mode is absent because the $\mathcal{M}_{N,k,j}^{({\boldsymbol{0}})}(\boldsymbol{\nu})$ term integrated to $0$ against the complex exponential over its
fundamental domain.

If we were to integrate positive Fourier modes over complete circles we would get $0$.
But here we sum over the discrete analogue of the circle: $\mathcal{V}_k$ and $\mathcal{V}_j$.
That means if we had a positive Fourier mode in the variable in $\mathcal{V}_k$ of a positive integer multiple of $k$,
or if we had a positive Fourier mode in the variable in $\mathcal{V}_j$ of a positive integer multiple of $j$, it would appear to be a
zero mode on the discrete lattice.
All modes beneath the $k$th or $j$th modes, for the respective variables ($y_1$ or $y_2$ in the former case and $y_3$ in the latter)
will sum to $0$ by the formula
$$
\sum_{n=0}^{R-1} e^{2\pi i n S}\, =\, \begin{cases} R & \text{ if $S \in R \Z = \{\dots,-R,0,R,\dots\}$,}\\ 0 & \text{ if $S \in \Z \setminus R\Z$}
\end{cases}
$$
Note that $k$ and $j$ are both of order $N$ since $k \sim \kappa N$ and $j \sim \gamma N$ as $N \to \infty$ (with $\gamma,\kappa>0$ fixed).
So we can apply the elementary Fourier analysis result that the $R$th Fourier coefficient decays as $O(1/R^n)$
for an $n$-times differentiable function.
Since the functions are analytic, we may apply this. Therefore, we get control of the small $\lambda$ domain once more.
This time we side-stepped the Euler-Maclaurin summation formula in multiple dimensions.

We will eschew further details in this document.
At a later date, we may return to consider more involved examples.
For example, we might consider the asymptotics beyond the large deviation rate function.
Then we would be motivated to include full details.
%
%

\baselineskip=12pt
\bibliographystyle{plain}

\end{document}